\documentclass[11pt]{article}%
\usepackage{amsfonts}
\usepackage{amsmath}
\usepackage{amssymb}
\usepackage{amsxtra}
\usepackage{graphicx}
\usepackage{geometry}
\usepackage{color}
\usepackage{colortbl}
\usepackage{caption}%
\providecommand{\U}[1]{\protect\rule{.1in}{.1in}}
\newtheorem{theorem}{Theorem}

\newtheorem{lemma}[theorem]{Lemma}

\newtheorem{proposition}[theorem]{Proposition}
\newtheorem{remark}[theorem]{Remark}

\numberwithin{equation}{section}
\ifx\pdfoutput\relax\let\pdfoutput=\undefined\fi
\newcount\msipdfoutput
\ifx\pdfoutput\undefined\else
\ifcase\pdfoutput\else
\ifx\paperwidth\undefined\else
\ifdim\paperheight=0pt\relax\else\pdfpageheight\paperheight\fi
\ifdim\paperwidth=0pt\relax\else\pdfpagewidth\paperwidth\fi
\fi\fi\fi
\begin{document}

\title{Scalar problems in junctions of rods and a plate. \\II. Self-adjoint extensions and simulation models.}
\author{R.Bunoiu\thanks{University of Lorraine, IECL, CNRS UMR 7502, 3 rue Augustin
Fresnel, 57073, Metz, France; email: renata.bunoiu@univ-lorraine.fr.},
G.Cardone\thanks{Universit\`{a} del Sannio, Department of Engineering, Corso
Garibaldi, 107, 82100 Benevento, Italy; email: giuseppe.cardone@unisannio.it},
S.A.Nazarov\thanks{Mathematics and Mechanics Faculty, St. Petersburg State
University, 198504, Universitetsky pr., 28, Stary Peterhof, Russia; Peter the
Great St. Petersburg State Polytechnical University, Polytechnicheskaya ul.,
29, St. Petersburg, 195251, Russia; Institute of Problems of Mechanical
Engineering RAS, V.O., Bolshoj pr., 61, St. Petersburg, 199178, Russia; email:
srgnazarov@yahoo.co.uk.}}
\maketitle

\begin{abstract}
\medskip In this work we deal with a scalar spectral mixed boundary value
problem in a spacial junction of thin rods and a plate. Constructing
asymptotics of the eigenvalues, we employ two equipollent asymptotic models
posed on the skeleton of the junction, that is, a hybrid domain. We, first,
use the technique of self-adjoint extensions and, second, we impose algebraic
conditions at the junction points in order to compile a problem in a function
space with detached asymptotics. The latter problem is involved into a
symmetric generalized Green formula and, therefore, admits the variational
formulation. In comparison with a primordial asymptotic procedure, these two
models provide much better proximity of the spectra of the problems in the
spacial junction and in its skeleton. However, they exhibit the negative
spectrum of finite multiplicity and for these "parasitic" eigenvalues we
derive asymptotic formulas to demonstrate that they do not belong to the
service area of the developed asymptotic models.

\medskip


\medskip


\end{abstract}

\section{Introduction\label{sect1}}

\subsection{Motivations\label{sect1.1}}

As it was observed in \cite{BuCaNa1}, an asymptotic expansion of a solution of
the Poisson scalar mixed boundary-value problem in a junction of thin rods and
a thin plate in a certain range of physical parameters gains the rational
dependence on the big parameter $|\ln h|$ where $h>0$ is a small parameter
characterizing the diameters of the rods and the thickness of the plate.
Similar asymptotic forms had been discovered for other elliptic problems
stated in domains with singular perturbations of the boundaries, see the books
\cite{MaNaPl, Ilin, KoMaMo}. The aim of this paper is to study the scalar
spectral problem associated to the Poisson problem studied in \cite{BuCaNa1}.

Based on the asymptotic procedure in \cite[\S 2.2.4, \S 5.5.2, \S 9.1.3]%
{MaNaPl} and \cite{na285, BuCaNa1}, it is quite predictable that the
asymptotic expansions of the eigenpairs "eigenvalue/eigenfunction" become much
more complicated than the one obtained for the solution of the Poisson problem
and purchase the holomorphic dependence on the parameter $|\ln h|^{-1}$. Such
sophistication of the asymptotic expansions and the lack of algorithms
allowing to clarify the appearing holomorphic functions make them almost
useless in applications, especially for combined numerical and asymptotic methods.

In \cite{Lions} J.-L. Lions announced as an open question an application of
the technique of self-adjoint extensions of differential operators for
modeling boundary-value problems with singular perturbations. This technique
has been employed in \cite{na344, na389, na576} and others to deal with
particular types of perturbations but at our knowledge, was never used before
for the study of spectral problems for junctions of thin domains with
different limit dimension like 1d and 2d for the rods and plate in our 3d junction.

\begin{figure}[ptb]
\begin{center}
\includegraphics[scale=0.55]{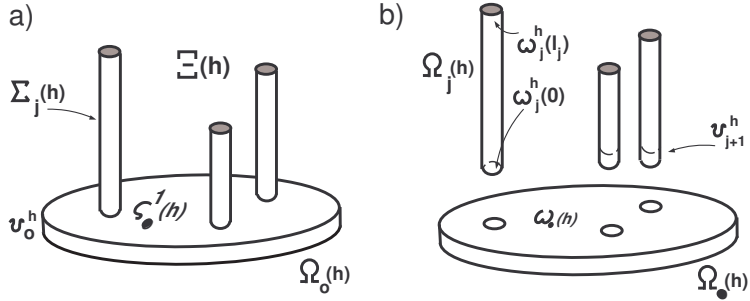}
\end{center}
\caption{The junction (a) and its elements (b).}%
\label{f1}%
\end{figure}

\begin{figure}[ptb]
\begin{center}
\includegraphics[scale=0.55]{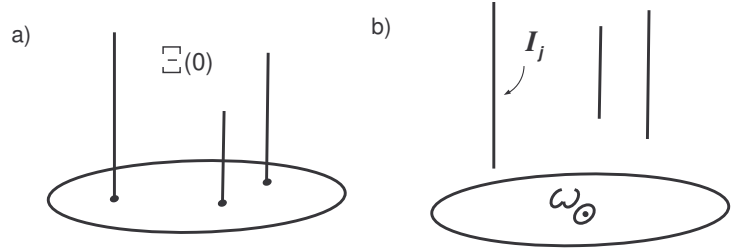}
\end{center}
\caption{The skeleton of the junction, the coupled (a) and disjoint (b) hybrid
domain.}%
\label{f2}%
\end{figure}

The use of such techniques in the asymptotic analysis here is the main novelty
of our work. We observe that to the spacial junction $\Xi(h)$ represented in
fig. \ref{f1},a, corresponds the hybrid domain $\Xi^{0}$ depicted in fig.
\ref{f2},a, which consists of several line segments joined to some interior
points of a planar domain. Supplying these elements of $\Xi^{0}$ with
differential structures, we describe a family of all self-adjoint operators
which is parametrized by a finite set of free parameters and choose an
appropriate set by examining the boundary layer phenomenon in the vicinity of
the junction zones, namely where the rods are inserted into small sockets in
the plate, fig. \ref{f1},b. As a result, we obtain a model which provides the
satisfactory proximity $O(h^{1/2}(1+|\ln h|)^{3})$ instead of the
uncomfortable one $O((1+|\ln h|)^{-2})$ within a simplified but comprehensible
version of the conventional asymptotic procedure. It should be mentioned that
our model involves only one scalar integral characteristic of each junction zone.

\subsection{Formulation of the spectral problem\label{sect1.2}}

Given a small parameter $h\in(0,h_{0}]$, we introduce the thin plate and rods%
\begin{align}
\Omega_{0}\left(  h\right)   &  =\left\{  x=\left(  y,z\right)  \in
\mathbb{R}^{3}:y=(y_{1},y_{2})\in\mathbb{\omega}_{0},\text{ }\zeta:=h^{-1}%
z\in\left(  0,1\right)  \right\}  ,\label{1}\\
\Omega_{j}\left(  h\right)   &  =\left\{  x:\eta^{j}=h^{-1}\left(
y-P^{j}\right)  \in\mathbb{\omega}_{j},\text{ \ }z\in I_{j}:=(0,l_{j}%
)\right\}  ,\ j=1,...,J. \label{2}%
\end{align}
Here, $\omega_{p},$ $p=0,1,...,J,$ are domains in the plane $\mathbb{R}^{2}$
with smooth (for simplicity) boundaries $\partial\mathbb{\omega}_{p}$ and the
compact closures $\overline{\mathbb{\omega}}_{p}=\mathbb{\omega}_{p}%
\cup\partial\mathbb{\omega}_{p};$ $l_{1},...,l_{J}$ are positive numbers
independent of $h$, and $P^{j}\in\omega_{0}$, $P^{j}\neq P^{k}$ for $j\neq k$.
Reducing a characteristic size of $\mathbb{\omega}_{0}$ to $1$, we make
Cartesian coordinates and all geometric parameters dimensionless. We fix some
$h_{0}\in(0,\min\left\{  l_{1},...,l_{J}\right\}  )$ such that, for
$h\in\left(  0,h_{0}\right]  ,$ $\overline{\mathbb{\omega}}_{j}^{h}%
\subset\mathbb{\omega}_{0}$ and $\overline{\mathbb{\omega}}_{j}^{h}%
\cap\overline{\mathbb{\omega}}_{k}^{h}=\emptyset,$ $j\neq k$, where
$\mathbb{\omega}_{j}^{h}=\left\{  y:\eta^{j}\in\mathbb{\omega}_{j}\right\}  $
is the cross section of the rod $\Omega_{j}\left(  h\right)  $ while
$\mathbb{\omega}_{j}^{h}(0)$ and $\mathbb{\omega}_{j}^{h}(l_{j})$ are its
lower and upper ends. In the sequel the bound $h_{0}$ can be diminished but
always remains strictly positive.

The rods (\ref{2}) are plugged into sockets, i.e., holes in the plate, fig.
\ref{f1},b and a,,%
\begin{equation}
\Omega_{\bullet}\left(  h\right)  =\mathbb{\omega}_{\bullet}\left(  h\right)
\times\left(  0,h\right)  ,\ \ \ \mathbb{\omega}_{\bullet}\left(  h\right)
=\mathbb{\omega}_{0}\mathbb{\diagdown}\left(  \overline{\mathbb{\omega}}%
_{1}^{h}\cup.....\cup\overline{\mathbb{\omega}}_{{J}}^{h}\right)  , \label{3}%
\end{equation}
and compose the junction%
\begin{equation}
\Xi\left(  h\right)  =\Omega_{\bullet}\left(  h\right)  \cup\Omega_{1}\left(
h\right)  \cup....\cup\Omega_{J}\left(  h\right)  . \label{4}%
\end{equation}
The lateral side of the plate (\ref{1}) and (\ref{3}) and their bases are
denoted by $\mathbb{\upsilon}_{0}\left(  h\right)  =\partial\omega_{0}%
\times(0,h)$ and%
\begin{equation}
\varsigma_{0}^{i}\left(  h\right)  =\left\{  x:y\in\mathbb{\omega}%
_{0},\ z=ih\right\}  ,\ \ \ \varsigma_{\bullet}^{i}\left(  h\right)  =\left\{
x:y\in\mathbb{\omega}_{\bullet},\ z=ih\right\}  ,\ \ i=0,1. \label{5}%
\end{equation}
The lateral side of the rod $\Omega_{j}\left(  h\right)  $ is divided into two
parts%
\[
\Sigma_{j}\left(  h\right)  =\partial\mathbb{\omega}_{j}^{h}\times\left(
h,l_{j}\right)  ,\ \ \mathbb{\upsilon}_{j}^{h}=\partial\mathbb{\omega}_{j}%
^{h}\times\left(  0,h\right)  ,
\]
where the latter is the junctional boundary of the socket.

In the junction (\ref{4}) we consider a spectral problem consisting of the
differential equations with the Laplace operator%
\begin{align}
-\Delta_{x}u_{0}\left(  h,x\right)   &  =\lambda(h)u_{0}\left(  h,x\right)
,\ \ \ x\in\Omega_{\bullet}\left(  h\right)  ,\label{7}\\
-\gamma_{j}\left(  h\right)  \Delta_{x}u_{j}\left(  h,x\right)   &
=\lambda(h)\rho_{j}(h)u_{j}\left(  h,x\right)  ,\ \ \ x\in\Omega_{j}\left(
h\right)  , \label{8}%
\end{align}
the Neumann and Dirichlet boundary conditions%
\begin{align}
\partial_{\nu}u_{0}\left(  h,x\right)   &  =0,\ \ \ \ x\in\Sigma_{\bullet
}\left(  h\right)  =\mathbb{\upsilon}_{0}\left(  h\right)  \cup\varsigma
_{\bullet}^{0}\left(  h\right)  \cup\varsigma_{\bullet}^{1}\left(  h\right)
,\label{9}\\
\gamma_{j}\left(  h\right)  \partial_{\nu}u_{j}\left(  h,x\right)   &
=0,\ \ \ x\in\Sigma_{j}\left(  h\right)  \cup\mathbb{\omega}_{j}^{h}\left(
0\right)  ,\label{10}\\
u_{j}\left(  h,x\right)   &  =0,\ \ \ x\in\mathbb{\omega}_{j}^{h}\left(
l_{j}\right)  , \label{11}%
\end{align}
and the transmission conditions%
\begin{align}
u_{0}\left(  h,x\right)   &  =u_{j}\left(  h,x\right)  ,\ \ \ \ x\in
\mathbb{\upsilon}_{j}^{h},\label{12}\\
\partial_{\nu}u_{0}\left(  h,x\right)   &  =\gamma_{j}\left(  h\right)
\partial_{\nu}u_{j}\left(  h,x\right)  ,\ \ \ x\in\mathbb{\upsilon}_{j}^{h}.
\label{13}%
\end{align}
Here, $\lambda(h)$ is a spectral parameter, $j=1,...,J,\ u_{0}$ and $u_{j}$
are restrictions of the function $u$ on $\Omega_{\bullet}\left(  h\right)  $
and $\Omega_{j}\left(  h\right)  $, respectively, and $\partial_{\nu}$ is the
outward normal derivative at the surface of $\overline{\Xi(h)}$ in (\ref{9})
and (\ref{10}) while in the second transmission condition $\partial_{\nu}$ is
outward with respect to the rods. In what follows we mainly deal with the
coefficients%
\begin{equation}
\gamma_{j}\left(  h\right)  =\gamma_{j}h^{-\alpha},\ \gamma_{j}%
>0,\ \ \ \ \ \rho_{j}\left(  h\right)  =\rho_{j}h^{-\alpha},\ \rho_{j}>0,
\label{14}%
\end{equation}
in the most informative but complicated case%
\begin{equation}
\alpha=1. \label{15}%
\end{equation}
We emphasize that in the case (\ref{15}) the limit passage $h\rightarrow+0$ in
the problem (\ref{7})-(\ref{13}) leads to the eigenvalue problem for a
self-adjoint operator in the skeleton of (\ref{4}) in fig. \ref{f2}a, a
hybrid domain, while, for $\alpha>1$ and $\alpha<1,$ this limit operator
decouples, cf fig. \ref{f2}b. We discuss the latter cases in Section
\ref{sect5} and pay attention to the homogeneous junction with
\begin{equation}
\alpha=0,\ \ \gamma_{j}=\rho_{j}=1,\ \ j=1,...,J. \label{82}%
\end{equation}
Other cases will be discussed in Section \ref{sect5}. The factors $\gamma_{j}$
and $\rho_{j}$ in (\ref{14}) are real positive numbers.

The variational formulation of the problem (\ref{7})-(\ref{13}) reads:%
\begin{equation}
a\left(  u,w;\Xi\left(  h\right)  \right)  =\lambda(h)b\left(  u,w;\Xi\left(
h\right)  \right)  \ \ \ \forall w\in H_{0}^{1}\left(  \Xi\left(  h\right)
;\Gamma\left(  h\right)  \right)  \label{16}%
\end{equation}
where $H_{0}^{1}\left(  \Xi\left(  h\right)  ;\Gamma\left(  h\right)  \right)
$ is the Sobolev space of functions satisfying the Dirichlet conditions
(\ref{11}) on $\Gamma\left(  h\right)  =\omega_{1}^{h}\left(  l_{1}\right)
\cup...\cup\omega_{J}^{h}\left(  l_{J}\right)  ,$%
\begin{align}
a\left(  u,v;\Xi\left(  h\right)  \right)   &  =\left(  \nabla_{x}u_{0}%
,\nabla_{x}w_{0}\right)  _{\Omega_{\bullet}\left(  h\right)  }+\sum
\nolimits_{j}\gamma_{j}\left(  h\right)  \left(  \nabla_{x}u_{j},\nabla
_{x}w_{j}\right)  _{\Omega_{j}\left(  h\right)  },\label{17}\\
b\left(  u,v;\Xi\left(  h\right)  \right)   &  =\left(  u_{0},w_{0}\right)
_{\Omega_{\bullet}\left(  h\right)  }+\sum\nolimits_{j}\rho_{j}\left(
h\right)  \left(  u_{j},w_{j}\right)  _{\Omega_{j}\left(  h\right)
},\nonumber
\end{align}
$\left(  \ ,\ \right)  _{\Xi\left(  h\right)  }$ is the natural scalar product
in the Lebesgue space $L^{2}\left(  \Xi\left(  h\right)  \right)  $ and
$\sum\nolimits_{j}$ everywhere stands for summation over $j=1,...,J$.

In view of the compact embedding $H^{1}\left(  \Xi\left(  h\right)  \right)
\subset L^{2}\left(  \Xi\left(  h\right)  \right)  $ the spectral problem
(\ref{16}) possesses, for every fixed $h$, the following positive monotone
unbounded sequence of eigenvalues
\begin{equation}
0<\lambda^{1}(h)<\lambda^{2}(h)\leq...\leq\lambda^{n}(h)\leq...\rightarrow
+\infty\label{18}%
\end{equation}
listed according to their multiplicity. The corresponding eigenfunctions
$u^{1}(h,\cdot),\ u^{2}(h,\cdot),...,\ u^{n}(h,\cdot),...\in H_{0}^{1}\left(
\Xi\left(  h\right)  ;\Gamma\left(  h\right)  \right)  $ can be subject to the
normalization and orthogonality conditions%
\begin{equation}
b\left(  u^{n},v^{m};\Xi\left(  h\right)  \right)  =\delta_{n,m}%
,\ \ n,m\in\mathbb{N}=\left\{  1,2,3,...\right\}  , \label{19}%
\end{equation}
where $\delta_{n,m}$ is the Kronecker symbol.

\subsection{The hybrid domain\label{sect1.3}}

The asymptotic analysis which has been presented at length, for example, in
\cite{BuCaNa1, ButCaNa1, ButCaNa2} and, in particular, includes the dimension
reduction procedure, converts the differential equations (\ref{7}) and
(\ref{8}) with the Neumann boundary conditions (\ref{9}) and (\ref{10}),
respectively, into the limit equations%
\begin{align}
-\Delta_{y}v_{0}(y)  &  =\mu v_{0}(y),\ \ \ y\in\omega_{\odot},\label{20}\\
-\gamma_{j}|\omega_{j}|\partial_{z}^{2}v_{j}(z)  &  =\mu\rho_{j}|\omega
_{j}|v_{j}(z),\ \ \ z\in(0,l_{j}), \label{21}%
\end{align}
where $\omega_{\odot}$ is the punctured domain $\omega_{0}\setminus
\mathcal{P},$ $\mathcal{P}=\left\{  P^{1},...,P^{J}\right\}  $ and
$|\omega_{j}|$ is the area of the domain $\omega_{j}$. Moreover, a primary
examination of the boundary layer phenomenon near the lateral side
$\mathbb{\upsilon}_{0}\left(  h\right)  $ of the plate and the end $\omega
_{j}^{h}(l_{j})$ of the rod $\Omega_{j}\left(  h\right)  $, respectively,
gives the following boundary conditions%
\begin{align}
\partial_{\nu}v_{0}(y)  &  =0,\ \ \ y\in\partial\omega_{0},\label{22}\\
v_{j}(l_{j})  &  =0. \label{23}%
\end{align}

However, the one-dimensional and two-dimensional problems are not completed
yet due to the lack of boundary conditions at the endpoints $z=0$ of the
intervals $(0,l_{j})$ and because the differential equation (\ref{20}) is
fulfilled for sure only outside the points $P^{1},...,P^{J}$, since near the
sockets $\omega_{j}^{h}\times(0,h)$ the geometrical structure of the junction
(\ref{4}) changes and becomes crucially spacial so that the dimension
reduction does not work. As was shown in \cite{BuCaNa1}, this observation
requires to consider solutions of the problem (\ref{20}), (\ref{22}) with
logarithmic singularities at the points in the set $\mathcal{P}$. We also will
take in the sequel such singular solutions into account, however further
considerations in this paper diverge from the asymptotic analysis used in
\cite{BuCaNa1}. Indeed, we will provide two abstract but applicable
formulations of a spectral problem in the hybrid domain in fig. \ref{f2},a,
which give an approximation of the spectrum (\ref{18}) with relatively high
precision. First, we detect a self-adjoint operator as an extension of the
differential operator of the problem (\ref{20})-(\ref{23}) supplied (cf.
(\ref{A3}) and (\ref{A1})) with the restrictive conditions%
\begin{equation}
v_{0}(P^{j})=0,\ \ v_{j}(0)=\partial_{z}v_{j}(0)=0,\ \ j=1,...,J. \label{too}%
\end{equation}

Second, we construct certain point conditions at $P^{1},...,P^{J}$ which tie
the independent problems (\ref{20}), (\ref{22}) and (\ref{21}), (\ref{23})
into a formally self-adjoint problem in the hybrid domain $\Xi^{0}.$ These two
formulations happen to be equivalent and both are realized as operators with
the discrete spectrum which, in the low-frequency range, approximate the
spectrum of the problem (\ref{7})-(\ref{13}) (or, equivalently, (\ref{16}))
with admissible precision\footnote{The error estimates are derived in the
paper with quite simple tools. Advanced estimation may detect the accuracy
$O(h|\ln h|)$ and extend the proximity property of the models to a part of the
mid-frequency range, cf. \cite{na576}. The latter, however, enlarges
enormously massif of calculations.} $O(h^{1/2}|\ln h|^{3})$. Unfortunately,
serving for a particular range of the spectrum, both the operators lose the
positivity property and gain so called \textit{parasite} eigenvalues which are
negative and big, of order $h^{-2}$; therefore, we prove that they lay outside
the scope of the asymptotic models and have no relation to the original
problem. In order to furnish, for example, an application of the minimum
principle, cf. \cite[Thm. 10.2.1]{BiSo}, we construct detailed asymptotics of
these \textit{parasite} eigenvalues and of the corresponding eigenfunctions,
which are located in the very vicinity of the points $P^{1},...,P^{J}$ and
decay exponentially at a distance from them.

Parameters of the self-adjoint extension and ingredients of the point
conditions are found out with the help of the method of matched asymptotic
expansions on the basis of special solutions described in the first part
\cite{BuCaNa1} of our work. Both linearly depend on $|\ln h|$ and this makes
the eigenvalues and eigenvectors to be real analytic functions in $|\ln
h|^{-1}.$ However, a possible numerical realization of the models with a small
but fixed parameter $h$ does not require to take into account such
complication of asymptotic expansions. We, of course, compute explicitly
couple of initial terms of the convergent series in $|\ln h|^{-1}.$

\subsection{Outline of the paper\label{sect1.4}}

In Section \ref{sect2} we describe all self-adjoint extensions of the operator
of the problem (\ref{20})-(\ref{too}) as well as the point conditions which
involve the problems (\ref{20}), (\ref{22}) and (\ref{21}), (\ref{23}) into a
symmetric generalized Green formula. This material is known and is presented
in a condensed form, mainly in order to introduce the notation and explain
some technicalities used throughout the paper. We refer to the review papers
\cite{Pav}, \cite{Pank} and \cite{na345} for a detailed information. If the
skeleton $\Xi(0)$ decouples in the limit, see fig. \ref{f2},b, then the
extended operator of the Neumann problem (\ref{20}), (\ref{22}) may require
for "potentials of zero radii" \cite{BeFa, Pav} or "pseudo-Laplacian" in the
terminology \cite{Yves}, but in this case the ordinary differential equations
(\ref{21}), (\ref{23}) are supplied with either Neumann, or Dirichlet
condition at the endpoints $z=0$ of the interval (cf., Remark \ref{remSIMPLE},
\ref{remDIRI} and see the paper \cite{Pank} which provides the complete
description of the techniques of self-adjoint extensions).

The most interesting situation occurs under the restriction (\ref{15}) when
the skeleton does not decouple in the limit, see fig. \ref{f2},a. The
asymptotic procedure in \cite{BuCaNa1} allow us to determine in Section
\ref{sect3} appropriate parameters of a particular self-adjoint extension
serving for the original problem (\ref{7})-(\ref{13}) as well as all
ingredients of the point condition in the corresponding hybrid model.

The most cumbersome part of our analysis is concentrated in Section
\ref{sect4}, where we perform the justification of our asymptotic models with
the help of weighted estimates obtained in \cite{BuCaNa1}. We state here the
main result of the paper, Theorem \ref{thASYM}. Section \ref{sect5} contains
some simple asymptotic formulas and the example of the homogeneous junction,
see (\ref{82}).

\section{General statement of problems on the hybrid domain\label{sect2}}

\subsection{Unbounded operators and their adjoints\label{sect2.1}}

Let $A_{j}$ be an unbounded operator in the Lebesgue space $L^{2}(I_{j})$ with
the differential expression $-\gamma_{j}|\omega_{j}|\partial_{z}^{2}$ and the
domain%
\begin{equation}
\mathcal{D}(A_{j})=\left\{  w_{j}\in H^{2}(I_{j}):v_{j}(l_{j})=0,\ v_{j}%
(0)=\partial_{z}v_{j}(0)=0\right\}  . \label{A1}%
\end{equation}
The operator is symmetric and closed. By a direct calculation, it follows that
the adjoint operator $A_{j}^{\ast}$ gets the same differential expression but
its domain is bigger, namely%
\begin{equation}
\mathcal{D}(A_{j}^{\ast})=\left\{  w_{j}\in H^{2}(I_{j}):v_{j}(l_{j}%
)=0\right\}  . \label{A2}%
\end{equation}
Hence, $\dim(\mathcal{D}(A_{j}^{\ast})/\mathcal{D}(A_{j}))=2$ and the defect
index of $A_{j}$ is $1:1$.

Analogously, we introduce the unbounded operator $A_{0}$ in $L^{2}(\omega
_{0})$ with the differential expression $-\Delta_{y}=-\partial^{2}/\partial
y_{1}^{2}-\partial^{2}/\partial y_{2}^{2}$ and the domain
\begin{equation}
\mathcal{D}(A_{0})=\left\{  w_{0}\in H^{2}(\omega_{0}):\partial_{\nu}%
w_{0}=0\text{ on }\partial\omega_{0},\ w_{j}(P^{k})=0,\ k=1,...,J\right\}  .
\label{A3}%
\end{equation}
By virtue of the Sobolev embedding theorem $H^{2}(\omega_{0})\subset
C(\omega_{0})$ and the classical Green formula%
\begin{equation}
-(\Delta_{y}w_{0},v_{0})_{\omega_{0}}+(\partial_{\nu}w_{0},v_{0}%
)_{\partial\omega_{0}}=-(w_{0},\Delta_{y}v_{0})_{\omega_{0}}+(w_{0}%
,\partial_{\nu}v_{0})_{\partial\omega_{0}}, \label{A4}%
\end{equation}
the operator $A_{0}$ is closed and symmetric. The following lemma, in
particular, shows that the defect index of $A_{0}$ is $J\times J$.

\begin{lemma}
\label{lem1A}The adjoint operator $A_{0}^{\ast}$ for $A_{0}$ has the
differential expression $-\Delta_{y}$ and the domain%
\begin{align}
\mathcal{D}(A_{0}^{\ast})  &  =\left\{  V_{0}\in L^{2}(\omega_{0}%
):V_{0}(y)=\widehat{V}_{0}(y)-\frac{1}{2\pi}%
%TCIMACRO{\tsum \nolimits_{j}}%
%BeginExpansion
{\textstyle\sum\nolimits_{j}}
%EndExpansion
b_{j}\chi_{j}(y)\ln|y-P^{j}|,\right.  \ \label{A5}\\
&  \ \ \ \ \ \ \ \ \ \ \ \ \ \ \ \ \ \ \ \ \ \ \ \ \ \ \ \ \ \ \ \left.
b_{j}\in\mathbb{C},\ \widehat{V}_{0}(y)\in H^{2}(\omega_{0}),\ \partial_{\nu
}\widehat{V}_{0}=0\text{ on }\partial\omega_{0}\right\}  ,\nonumber
\end{align}
where $\chi_{1},...,\chi_{J}\in C_{c}^{\infty}(\omega_{0})$ are cut-off
functions such that
\[
\chi_{j}(P^{j})=1,\ \ \ \chi_{j}(y)\chi_{k}(y)=0\text{ for }j\neq
k,\ \ \ \mathrm{supp}\chi_{j}\subset\omega_{0}.
\]

\end{lemma}

\textbf{Proof.} By definition, a function $V_{0}\in L^{2}(\omega_{0})$ belongs
to $\mathcal{D}(A_{0}^{\ast})$ if and only if the following integral identity
holds:%
\begin{equation}
-(V_{0},\Delta_{y}v_{0})_{\omega_{0}}=(F_{0},v_{0})_{\omega_{0}}\ \ \ \forall
v_{0}\in\mathcal{D}(A_{0}). \label{A6}%
\end{equation}
At the first step we take $v_{0}\in C_{c}^{\infty}(\overline{\omega}%
_{0}\setminus\mathcal{P})\cap\mathcal{D}(A_{0}).$ Based on the Green formula
(\ref{A4}), we recall the Neumann boundary condition in (\ref{A3}) and apply
classical results \cite[\S 3-6, Ch. 2]{LiMa} on lifting smoothness of
solutions to elliptic problems. In this way we conclude that $V_{0}\in
H_{loc}^{2}(\overline{\omega}_{0}\setminus\mathcal{P})$ and%
\begin{equation}
-\Delta_{y}V_{0}(y)=F_{0}(y),\ y\in\omega_{\odot}=\omega_{0}\setminus
\mathcal{P},\ \ \ \ \partial_{\nu}V_{0}(y)=0\text{, }y\in\partial\omega_{0}.
\label{A7}%
\end{equation}

The next step requires for results \cite{Ko} of the theory of elliptic problem
in domains with conical points. Indeed, regarding $P^{j}$ as the top of the
"complete cone" $\mathbb{R}^{2}\setminus P^{j}$, that is, the punctured plane,
we introduce the Kondratiev space $V_{\beta}^{l}(\omega_{0})$ with
$l\in\left\{  0,1,2,...\right\}  $ and $\beta\in\mathbb{R}$ as the completion
of $C_{c}^{\infty}(\overline{\omega}_{0}\setminus\mathcal{P})$ with
respect to the weighted norm%
\begin{equation}
||v_{0};V_{\beta}^{l}(\omega_{0})||=\left(
%TCIMACRO{\tsum \limits_{k=0}^{l}}%
%BeginExpansion
{\textstyle\sum\limits_{k=0}^{l}}
%EndExpansion
||\min\left\{  r_{1},...,r_{J}\right\}  ^{\beta-l-k}\nabla_{y}^{k}v_{0}%
;L^{2}(\omega_{0})^{2}||\right)  ^{1/2} \label{A8}%
\end{equation}
where $r_{j}=|y-P^{j}|$ and $\nabla_{y}^{k}v_{0}$ stands for a collection of
all order $k$ derivatives of $v_{0}$. Clearly, $L^{2}(\omega_{0})\subset
V_{\delta}^{0}(\omega_{0})$ and $V_{\delta}^{2}(\omega_{0})\subset
H^{1}(\omega_{0})$ for $\delta\in\lbrack0,1]$.

Since $V_{0}\in L^{2}(\omega_{0})\subset V_{0}^{0}(\omega_{0})$, the theorem
on asymptotics \cite{Ko}, see, e.g., \cite[\S 3.5, \S 4.2, \S 6.4]{NaPl} and
also the introductory chapters in the books \cite{NaPl, KoMaRo1}, gives the
representation%
\begin{equation}
V_{0}(y)=\widetilde{V}_{0}(y)+%
%TCIMACRO{\tsum \nolimits_{j}}%
%BeginExpansion
{\textstyle\sum\nolimits_{j}}
%EndExpansion
\chi_{j}(y)\left(  a_{j}-\frac{b_{j}}{2\pi}\ln r_{j}\right)  \label{A9}%
\end{equation}
as well as the inclusion $\widetilde{V}_{0}\in V_{\delta}^{2}(\omega_{0})$
with any $\delta>0$ and the estimate%
\begin{equation}
||\widetilde{V}_{0};V_{\delta}^{2}(\omega_{0})||+%
%TCIMACRO{\tsum \nolimits_{j}}%
%BeginExpansion
{\textstyle\sum\nolimits_{j}}
%EndExpansion
(|a_{j}|+|b_{j}|)\leq c_{\delta}(\left\Vert F_{0};L^{2}(\omega_{0})\right\Vert
+\left\Vert V_{0};L^{2}(\omega_{0})\right\Vert ). \label{A10}%
\end{equation}
Hence, the sum%
\begin{equation}
\widehat{V}_{0}(y)=V_{0}(y)+\frac{1}{2\pi}\sum\nolimits_{j}b_{j}\ln r_{j}
\label{A99}%
\end{equation}
belongs to $H^{1}(\omega_{0})$ and still solves the problem (\ref{A7}) with a
new right-hand side $\widehat{F}_{0}\in L^{2}(\omega_{0})$ in the Poisson
equation. Thus, referring to \cite[\S 9, Ch. 2]{LiMa} we have $\widehat{V}%
_{0}\in H^{2}(\omega_{0})$ and, therefore, $V_{0}$ falls into the linear set
(\ref{A5}). $\blacksquare$

\begin{remark}
\label{remFOURIER}The representation (\ref{A9}) can be derived by means of the
Fourier method but the Kondratiev theory \cite{Ko} helps to avoid any calculation.
\end{remark}

\subsection{The generalized Green formula\label{sect2.2}}

The norm, cf. the left-hand side of (\ref{A10}),%
\begin{equation}
\left\Vert V_{0};\mathfrak{H}_{0}\right\Vert =(||\widehat{V}_{0};H^{2}%
(\omega_{0})||^{2}+%
%TCIMACRO{\tsum \nolimits_{j}}%
%BeginExpansion
{\textstyle\sum\nolimits_{j}}
%EndExpansion
|b_{j}|^{2})^{1/2} \label{A11}%
\end{equation}
brings Hilbert structure to the linear space (\ref{A5}) denoted by
$\mathfrak{H}_{0}$. By $\mathfrak{H}$, we understand the direct product%
\begin{equation}
\mathfrak{H}=\mathfrak{H}_{0}\times\mathfrak{H}_{1}\times...\times
\mathfrak{H}_{J} \label{A12}%
\end{equation}
where $\mathfrak{H}_{j}$ is the linear space (\ref{A2}) with the Sobolev
$H^{2}$-norm. Moreover, taking into account the right-hand sides of the
equations (\ref{20}) and (\ref{21}) we supply the vector Lebesgue space
$\mathfrak{L}$ with the special norm%
\begin{equation}
||v;\mathfrak{L}||=(||v_{0};L^{2}(\omega_{0})||^{2}+%
%TCIMACRO{\tsum \nolimits_{j}}%
%BeginExpansion
{\textstyle\sum\nolimits_{j}}
%EndExpansion
\rho_{j}|\omega_{j}|\ ||v_{j};L^{2}(I_{j})||^{2})^{1/2}, \label{LL}%
\end{equation}
with $v=(v_{0},v_{1},...,v_{J})\in\mathfrak{L}:=L^{2}(\omega_{0})\times
L^{2}(I_{1})\times...\times L^{2}(I_{J}).$

We also introduce two continuous projections $\wp_{\pm}:\mathfrak{H}%
\rightarrow\mathbb{R}^{2J}$ by the formulas%
\begin{align}
\wp_{+}v  &  =(\wp_{+}^{\prime}v,\wp_{+}^{\prime\prime}v)=(\widehat{v}%
_{0}(P^{1}),...,\widehat{v}_{0}(P^{J}),v_{1}(0),...,v_{J}(0)),\label{A13}\\
\wp_{-}v  &  =(\wp_{-}^{\prime}v,\wp_{-}^{\prime\prime}v)=(b_{1}%
,...,b_{J},-\gamma_{1}|\omega_{1}|\partial_{z}v_{1}(0),...,-\gamma_{J}%
|\omega_{J}|\partial_{z}v_{J}(0)),\nonumber
\end{align}
where $v=(v_{0},v_{1},...,v_{J})\in\mathfrak{H}$ and $b_{1},...,b_{J}$,
$\widehat{v}_{0}$ are attributes of the decomposition (\ref{A9}) of $v_{0}%
\in\mathfrak{H}_{0}$.

The next assertion is but a concretization of a general result in \cite{na161,
na159, na165}, see also \cite[\S 6.2]{NaPl}, however we give a condensed and
much simplified proof for reader's convenience.

\begin{proposition}
\label{prop2A}For $v$ and $w=(w_{0},w_{1},...,w_{J})$ in $\mathfrak{H}$, the
generalized Green formula%
\begin{align}
q(v,w)  &  :=-(\Delta_{y}v_{0},w_{0})_{\omega_{0}}+(v_{0},\Delta_{y}%
w_{0})_{\omega_{0}}-%
%TCIMACRO{\tsum \nolimits_{j}}%
%BeginExpansion
{\textstyle\sum\nolimits_{j}}
%EndExpansion
\gamma_{j}|\omega_{j}|((\partial_{z}^{2}v_{j},w_{j})_{I_{j}}-(v_{j}%
,\partial_{z}^{2}w_{j})_{I_{j}})\label{A14}\\
&  =\left\langle \wp_{+}v,\wp_{-}w\right\rangle -\left\langle \wp_{-}v,\wp
_{+}w\right\rangle \nonumber
\end{align}
is valid, where $\left\langle \ ,\ \right\rangle $ stands for the natural
scalar product in $\mathbb{R}^{2J}$.
\end{proposition}

\textbf{Proof.} First of all, we write the evident identity%
\begin{equation}
-(\partial_{z}^{2}v_{j},w_{j})_{I_{j}}+(v_{j},\partial_{z}^{2}w_{j})_{I_{j}%
}=w_{j}(0)\partial_{z}v_{j}(0)-v_{j}(0)\partial_{z}w_{j}(0) \label{A15}%
\end{equation}
and multiply it with $\gamma_{j}|\omega_{j}|$. Then we take $v_{0}%
,w_{0}\in\mathfrak{H}_{0}$ and write the standard Green formula
\begin{align}
-(  &  \Delta_{y}v_{0},w_{0})_{\omega_{0}}+(v_{0},\Delta_{y}w_{0})_{\omega
_{0}}=\lim_{\delta\rightarrow+0}(-(\Delta_{y}v_{0},w_{0})_{\omega_{\delta}%
}+(v_{0},\Delta_{y}w_{0})_{\omega_{\delta}})\label{A16}\\
&  =-\lim_{\delta\rightarrow+0}\sum\nolimits_{j}\delta\int_{0}^{2\pi}\left(
\left(  \widehat{w}(P^{j})-\frac{b_{j}^{w}}{2\pi}\ln r_{j}\right)
\frac{\partial}{\partial r_{j}}\frac{b_{j}^{v}}{2\pi}\ln r_{j}-\left(
\widehat{v}(P^{j})-\frac{b_{j}^{v}}{2\pi}\ln r_{j}\right)  \frac{\partial
}{\partial r_{j}}\frac{b_{j}^{w}}{2\pi}\ln r_{j}\right)  d\varphi
_{j}\nonumber\\
&  =-\sum\nolimits_{j}(\widehat{w}(P^{j})b_{j}^{v}-\widehat{v}(P^{j})b_{j}%
^{w})\nonumber
\end{align}
where $\omega_{\delta}=\omega_{0}\setminus\mathbb{B}_{\delta}(P^{j}),$
$\mathbb{B}_{\delta}(P^{j})=\{y:r_{j}<\delta\}$ is a disk and $(r_{j}%
,\varphi_{j})\in\mathbb{R}_{+}\times\lbrack0,2\pi)$ is the polar coordinate
system centered at $P^{j}$. Now (\ref{A14}) follows from (\ref{A15}),
(\ref{A16}) and (\ref{A13}). $\blacksquare$

\subsection{Self-adjoint extensions\label{sect2.3}}

Calculations in Section \ref{sect2.2} detect the defect index $2J:2J$ of the
operator $A=(A_{0},A_{1},...,A_{J})$ with the differential expression
$(-\Delta_{y}v_{0},-\gamma_{1}\rho_{1}^{-1}\partial_{z}^{2},...,-\gamma
_{J}\rho_{J}^{-1}\partial_{z}^{2})$ in the Hilbert space $\mathfrak{L}$, see
(\ref{LL}). Hence, this operator admits a self-adjoint extension $\mathcal{A}%
$, that is, $A\subset\mathcal{A}\subset A^{\ast}$ and $\mathcal{A}%
=\mathcal{A}^{\ast}$.

Since $\mathcal{A}$ is a restriction of $A^{\ast}=(A_{0}^{\ast},A_{1}^{\ast
},...,A_{J}^{\ast})$, we conclude that%
\begin{equation}
\mathcal{D}(\mathcal{A})=\mathcal{D}(\mathcal{A}^{\ast})=\left\{
v\in\mathcal{D}(\mathcal{A}):(\wp_{+}v,\wp_{-}v)\in\mathcal{R}\right\}
\label{A17}%
\end{equation}
where a linear subspace $\mathcal{R}\subset\mathbb{R}^{4J}$ of dimension $2J$
must be chosen such that in accord with Proposition \ref{prop2A}%
\begin{equation}
q(v,w)=0\ \ \ \forall v,w\in\mathcal{D}(\mathcal{A}). \label{A18}%
\end{equation}

The symplectic form (\ref{A14}) is actually defined on the factor space
$\mathcal{D}(\mathcal{A}^{\ast})/\mathcal{D}(\mathcal{A})\approx
\mathbb{R}^{4J}$ because the generalized Green formula demonstrates that%
\begin{equation}
0=q(v,w)=-\overline{q(w,v)}\text{ \ for }v\in\mathcal{D}(\mathcal{A}%
),\ w\in\mathcal{D}(\mathcal{A}^{\ast}). \label{A19}%
\end{equation}
Description of null spaces of a symplectic form in Euclidean space is a
primary algebraic question, cf. \cite{Leng}, but it gives a direct
identification of all self-adjoint extensions of our operator $A$, see
\cite{BeFa}, \cite{Rofe}, \cite{Yves} and, e.g., \cite{Pav, Pank, na239}.

\begin{proposition}
\label{prop3A}Let $\mathcal{R}^{+}\mathcal{\oplus R}^{0}\mathcal{\oplus R}%
^{-}$ be an orthogonal decomposition of $\mathbb{R}^{2J}$ and let
$\mathcal{S}$ be a symmetric invertible operator in $\mathcal{R}^{0}$. The
restriction $\mathcal{A}$ of the operator $A^{\ast}$ onto the domain%
\begin{equation}
\mathcal{D}(\mathcal{A})=\left\{  v\in\mathcal{D}(A^{\ast}):\wp_{+}%
v=t^{+}+\mathcal{T}t^{0},\ \wp_{-}v=t^{-}+t^{0},\ t^{\alpha}\in\mathcal{R}%
^{\alpha},\ \alpha=0,\pm\right\}  \label{A20}%
\end{equation}
is a self-adjoint extension of the operator $A$ in $\mathfrak{L}$. Any
self-adjoint extension of $A$ can be obtained in this way.
\end{proposition}

\begin{remark}
\label{remSIMPLE}If we put
\begin{equation}
\mathcal{R}^{-}=\mathcal{R}^{0}=\{0\}^{2J},\ \ \mathcal{R}^{+}=\mathbb{R}%
^{2J}, \label{A21}%
\end{equation}
then the self-adjoint extension $\mathcal{A}^{0}$ given in Proposition
\ref{prop3A} is nothing but the set of the two-dimensional Neumann problem%
\begin{equation}
-\Delta_{y}v_{0}(y)=f_{0}(y),\ \ y\in\omega_{0},\ \ \ \ \partial_{\nu}%
v_{0}(y)=0,\ \ y\in\partial\omega_{0}, \label{A22}%
\end{equation}
and the one-dimensional mixed boundary-value problems%
\begin{equation}
-\gamma_{j}|\omega_{j}|\partial_{z}^{2}v_{j}(z)=f_{j}(z),\ z\in(0,l_{j}%
),\ \ \ v_{j}(l_{j})=\gamma_{j}|\omega_{j}|\partial_{z}v_{j}(0)=0. \label{A23}%
\end{equation}
These problems are independent and are posed on the spaces $H^{2}(\omega_{0})$
and $H^{2}(I_{j})$, respectively. The corresponding operator has the kernel
spanned over the constant vectors $(c_{0},0,...,0)$. $\blacksquare$
\end{remark}

\begin{remark}
\label{remDIRI}In the case $\mathcal{R}^{-}=\{0\}^{J}\times\mathbb{R}^{J},$
$\mathcal{R}^{0}=\{0\}^{2J},\ \mathcal{R}^{+}=\mathbb{R}^{J}\times\{0\},$ we
obtain a self-adjoint extension which gives rise to the Neumann problem
(\ref{A22}) and a set of the Dirichlet problems%
\[
-\gamma_{j}|\omega_{j}|\partial_{z}^{2}v_{j}(z)=f_{j}(z),\ z\in(0,l_{j}%
),\ \ \ v_{j}(l_{j})=v_{j}(0)=0.\ \ \blacksquare
\]

\end{remark}

In Section \ref{sect3} we come across a self-adjoint extension $\mathcal{A}%
=\mathcal{A}^{h}$ (the superscript will appear in Section \ref{sect3.2}) with
the following attributes in (\ref{A20}):%
\begin{equation}
\mathcal{R}^{-}=\{0\}^{2J},\ \mathcal{R}^{0}=\left\{  \binom{c}{-c}%
:c\in\mathbb{R}^{J}\right\}  ,\ \mathcal{R}^{+}=\left\{  \binom{c}{c}%
:c\in\mathbb{R}^{J}\right\}  ,\ \mathcal{T}=\left(
\begin{array}
[c]{cc}%
\mathbb{O} & \frac{1}{2}S\\
\frac{1}{2}S & \mathbb{O}%
\end{array}
\right)  \label{A24}%
\end{equation}
where $S=S^{h}$ is a symmetric non-degenerate matrix of size $J\times J$. In
the next section we will give a different formulation of the abstract equation%
\begin{equation}
\mathcal{A}^{h}v=f\in\mathcal{L} \label{A25}%
\end{equation}
which will help us to study the spectrum of $\mathcal{A}^{h}.$

\subsection{The differential problem with point conditions\label{sect2.4}}

Let $\wp_{\pm}^{\prime}$ and $\wp_{\pm}^{\prime\prime}$ be projections
$:\mathfrak{H}\rightarrow\mathbb{R}^{J}$ defined in (\ref{A13}). Following
\cite{na188, na239, na504}, we rewrite relations imposed on $\wp_{\pm}%
^{\prime}v$ and $\wp_{\pm}^{\prime\prime}v$ in (\ref{A20}) according to
(\ref{A24}), as the \textit{point conditions}
\begin{align}
\wp_{+}^{\prime\prime}v-\wp_{+}^{\prime}v-S\wp_{-}^{\prime}v  &
=0\in\mathbb{R}^{J},\label{A26}\\
\wp_{-}^{\prime}v+\wp_{-}^{\prime\prime}v  &  =0\in\mathbb{R}^{J}. \label{A27}%
\end{align}
We also will deal with the inhomogeneous equation%
\begin{equation}
\wp_{+}^{\prime\prime}v-\wp_{+}^{\prime}v-S\wp_{-}^{\prime}v=k\in
\mathbb{R}^{J}. \label{A28}%
\end{equation}
The problems%
\begin{align}
-\Delta_{y}v_{0}(y)  &  =f_{0}(y),\ \ y\in\omega_{\odot},\ \ \ \ \partial
_{\nu}v_{0}(y)=0,\ \ y\in\partial\omega\label{A29}\\
-\gamma_{j}|\omega_{j}|\partial_{z}^{2}v_{j}(z)  &  =f_{j}(z),\ z\in
(0,l_{j}),\ \ \ v_{j}(l_{j})=0 \label{A30}%
\end{align}
with the point conditions (\ref{A27}), (\ref{A28}) give rise to the continuous
mapping%
\begin{equation}
\mathfrak{A}:\mathfrak{H}_{-}=\{v\in\mathfrak{H}:\wp_{-}^{\prime}v+\wp
_{-}^{\prime\prime}v=0\}\rightarrow\mathfrak{R}:=\mathfrak{L}\times
\mathbb{R}^{J}. \label{A31}%
\end{equation}

\begin{remark}
\label{remEQU}Simple algebraic transformations demonstrate that, under
circumstances (\ref{A20}) and (\ref{A24}), a solution of the equation
(\ref{A25}) is a solution of the problem (\ref{A29}), (\ref{A30}),
(\ref{A26}), (\ref{A27}) and vice versa.
\end{remark}

\begin{proposition}
\label{prop5A}The operator $\mathfrak{A}$ in (\ref{A31}) is Fredholm of index zero.
\end{proposition}

\textbf{Proof. }The point conditions $\wp_{-}^{\prime}v=0,$ $\wp_{+}%
^{\prime\prime}v=0$ in Remark \ref{remDIRI} generate the Fredholm operator of
index zero%
\begin{equation}
H^{2}(\omega_{0})\times%
%TCIMACRO{\tprod \nolimits_{j=1}^{J}}%
%BeginExpansion
{\textstyle\prod\nolimits_{j=1}^{J}}
%EndExpansion
(H^{2}(I_{j})\cap H_{0}^{1}(I_{j}))\rightarrow\mathfrak{L}. \label{A32}%
\end{equation}
The operator (\ref{A31}) is a finite dimensional, i.e. compact, perturbation
of (\ref{A32}) and thus keeps the Fredholm property. Since $S=S^{T}$, the
generalized Green formula (\ref{A14}) can be written in the symmetric form
reflecting the particular point conditions (\ref{A26}), (\ref{A27})%
\begin{align}
q(v,w)  &  =\left\langle \wp_{+}^{\prime}v-\wp_{+}^{\prime\prime}v+S\wp
_{-}^{\prime}v,\wp_{-}^{\prime}w\right\rangle -\left\langle \wp_{-}^{\prime
}v,\wp_{+}^{\prime}w-\wp_{+}^{\prime\prime}w+S\wp_{-}^{\prime}w\right\rangle
\label{A33}\\
&  +\left\langle \wp_{+}^{\prime\prime}v,\wp_{-}^{\prime}w+\wp_{-}%
^{\prime\prime}w\right\rangle -\left\langle \wp_{-}^{\prime}v+\wp_{-}%
^{\prime\prime}v,\wp_{+}^{\prime\prime}w\right\rangle \nonumber
\end{align}
and hence an argument in \cite[Sect. 2.2.5, 2.5.3]{LiMa}, cf. \cite[\S 6.2]%
{NaPl}, shows that%
\begin{equation}
\mathrm{Ind~}\mathfrak{A}=\dim\ker\mathfrak{A}-\dim\text{\textrm{coker}%
}~\mathfrak{A}=0,\ \ \text{\textrm{coker}}~\mathfrak{A}=\{(v,\wp_{-}^{\prime
}v)\in\mathfrak{R}:v\in\ker\mathfrak{A}\}.\ \ \blacksquare\label{A333}%
\end{equation}

Let $G$ be the generalized Green function of the Neumann problem (\ref{A22}),
see, e.g., \cite{Smir}, namely a distributional solution of%
\begin{gather}
-\Delta_{y}G(y,\mathbf{y})=\delta(y-\mathbf{y})-|\omega_{0}|^{-1},\ y\in
\omega_{0},\ \ \ \partial_{\nu(y)}G(y,\mathbf{y})=0,\ y\in\partial\omega
_{0},\label{A34}\\
\int_{\omega_{0}}G(y,\mathbf{y})dy=0,\ \mathbf{y}\in\omega_{0},\nonumber
\end{gather}
where $\delta$ is the Dirac mass. We put $G^{j}(y)=G(y,P^{j})$ and write%
\begin{equation}
G^{j}(y)=-\chi_{j}(y)(2\pi)^{-1}\ln r_{j}+\widehat{G}^{j}(y),\ \ \ \widehat
{G}^{j}\in H^{2}(\omega_{0}). \label{A35}%
\end{equation}
The $J\times J$-matrix $\mathcal{G}$ with entries $\mathcal{G}_{k}%
^{j}=\widehat{G}^{j}(P^{k})$ is symmetric, see \cite[ Sect. 2.2]{BuCaNa1}. We
compose the vectors%
\begin{equation}
\mathbf{G}^{j}=(G^{j},\delta_{j1}\gamma_{1}^{-1}|\omega_{1}|^{-1}%
(z-l_{1}),...,\delta_{jJ}\gamma_{J}^{-1}|\omega_{J}|^{-1}(z-l_{J}%
))\in\mathfrak{H}, \label{A36}%
\end{equation}
which fall into the subspace $\mathfrak{H}_{-},$ see (\ref{A31}), because%
\begin{align}
\wp_{-}^{\prime}\mathbf{G}^{j}  &  =-\wp_{-}^{\prime\prime}\mathbf{G}%
^{j}=\mathbf{e}_{(j)},\ \label{A37}\\
\wp_{+}^{\prime}\mathbf{G}^{j}  &  =\mathcal{G}^{j}=(\mathcal{G}_{1}%
^{j},...,\mathcal{G}_{J}^{j}),\ \ \ \wp_{+}^{\prime\prime}\mathbf{G}%
^{j}=-\gamma_{j}^{-1}|\omega_{j}|^{-1}l_{j}\mathbf{e}_{(j)}.\nonumber
\end{align}
Here, $\mathbf{e}_{(j)}=(\delta_{j1},...,\delta_{jJ}),$ $j=1,...,J,$ is the
natural basis in $\mathbb{R}^{J}$.

Let $\mathcal{E}$ be a subspace spanned over the vector $\varepsilon
=(1,...,1)\in\mathbb{R}^{J},$ $|\varepsilon|=\sqrt{J},$ and $\mathbb{R}%
^{J}=\mathcal{E\oplus E}^{\bot}$ with the orthogonal projector $\mathcal{P}%
^{\bot}$ onto $\mathcal{E}^{\bot},$ $\dim\mathcal{E}^{\bot}=J-1$. We also
introduce the diagonal matrix%
\begin{equation}
\mathcal{Q}=\mathrm{diag}\{\gamma_{1}|\omega_{1}|l_{1}^{-1},...,\gamma
_{J}|\omega_{J}|l_{J}^{-1}\}. \label{A38}%
\end{equation}

\begin{theorem}
\label{th6A}If the operator
\begin{equation}
\mathcal{P}^{\bot}(S+\mathcal{G}+\mathcal{Q}^{-1})\mathcal{P}^{\bot
}:\mathcal{E}^{\bot}\rightarrow\mathcal{E}^{\bot} \label{AAA}%
\end{equation}
is invertible, then the problem (\ref{A27})-(\ref{A30}) on the hybrid domain
$\Xi^{0}=\omega_{\odot}\cup%
%TCIMACRO{\tbigcup \nolimits_{j=1}^{J}}%
%BeginExpansion
{\textstyle\bigcup\nolimits_{j=1}^{J}}
%EndExpansion
(P^{j}\cup I_{j})$, fig. \ref{f2},a, has a unique solution $v\in
\mathfrak{H}_{-}$ for any $\{f,k\}\in\mathfrak{R}$. In other words, the
operator (\ref{A31}) is an isomorphism.
\end{theorem}

\textbf{Proof. }We search for a solution of the problem in the form%
\begin{equation}
v=\mathbf{v}+\alpha_{1}\mathbf{G}^{1}+...+\alpha_{J}\mathbf{G}^{J} \label{A39}%
\end{equation}
where $\alpha=(\alpha_{1},...,\alpha_{J})\in\mathbb{R}^{J},$ $\mathbf{G}^{j}$
is given in (\ref{A36}) and $\mathbf{v}=(\mathbf{v}_{0},\mathbf{v}%
_{1},...,\mathbf{v}_{J})\in\mathfrak{H}_{-}$ with%
\begin{equation}
\mathbf{v}_{0}(y)=\alpha_{0}^{0}+\mathbf{v}_{0}^{0}(y),\ \mathbf{v}_{0}^{0}\in
H^{2}(\omega_{0}),\ \int_{\omega_{0}}\mathbf{v}_{0}^{0}(y)dy=0,\ \mathbf{v}%
_{j}\in H^{2}(I_{j}). \label{A40}%
\end{equation}
In view of (\ref{A34})-(\ref{A36}) these functions must satisfy the problems%
\begin{align}
-\Delta_{y}\mathbf{v}_{0}(y)  &  =f_{0}(y)-|\omega_{0}|^{-1}(\alpha
_{1}+...+\alpha_{J}),\ \ y\in\omega_{0},\ \ \ \ \partial_{\nu}\mathbf{v}%
_{0}(y)=0,\ \ y\in\omega_{0},\label{A41}\\
-\gamma_{j}|\omega_{j}|\partial_{z}^{2}\mathbf{v}_{j}(z)  &  =f_{j}%
(z),\ z\in(0,l_{j}),\ \ \ \ \mathbf{v}_{j}(l_{j})=0,\ \ \ -\gamma_{j}%
|\omega_{j}|\partial_{z}\mathbf{v}_{j}(0)=0.\nonumber
\end{align}
Under the condition%

\begin{equation}%
%TCIMACRO{\tsum \nolimits_{j}}%
%BeginExpansion
{\textstyle\sum\nolimits_{j}}
%EndExpansion
\alpha_{j}=|\omega_{0}|\int_{\omega_{0}}f_{0}(y)dy, \label{A42}%
\end{equation}
the problems (\ref{A41}) have unique solutions (\ref{A40}) but with arbitrary
constant $\alpha_{0}^{0}$. The vector function (\ref{A39}) fulfils the point
condition (\ref{A27}) while (\ref{A28}) turns into%
\begin{equation}
-\mathcal{G}\alpha-\alpha_{0}^{0}\varepsilon-\mathcal{Q}^{-1}\alpha
-S\alpha=k+\wp_{+}^{\prime}(\mathbf{v}_{0}^{0},0,...,0)-\wp_{+}^{\prime\prime
}(0,\mathbf{v}_{1},...,\mathbf{v}_{J})\in\mathbb{R}^{J}. \label{A43}%
\end{equation}
Applying the projector $\mathcal{P}^{\bot}$, we annul the term $\alpha_{0}%
^{0}\varepsilon$ in (\ref{A43}) and determine $\mathcal{P}^{\bot}\alpha$,
thanks to our assumption on the mapping (\ref{AAA}). Then the equation
(\ref{A42}) gives the remaining part of the coefficient vector in (\ref{A39}).
Recalling (\ref{A43}) yields a value of $\alpha_{0}^{0}$.

Since we have found a solution (\ref{A39}), the operator (\ref{A31}) is an
epimorphism and becomes isomorphism by virtue of Proposition \ref{prop5A}.
$\blacksquare$

\subsection{The variational formulation of the problem with point
conditions\label{sect2.5}}

Similarly to \cite{na188, na239, na504} we associate the problem
(\ref{A27})-(\ref{A30}) with the quadratic form%
\begin{align}
\mathfrak{E}(v;f,k)  &  =-\frac{1}{2}(\Delta_{y}v_{0},v_{0})_{\omega_{0}%
}-(f_{0},v_{0})_{\omega_{0}}-\sum\nolimits_{j}(\gamma_{j}|\omega_{j}%
|(\partial_{z}^{2}v_{j},v_{j})_{I_{j}}+(f_{j},v_{j})_{I_{j}}\label{A44}\\
&  +\frac{1}{2}\left\langle \wp_{+}^{\prime\prime}v-\wp_{+}^{\prime}v-S\wp
_{-}^{\prime}v,\wp_{-}^{\prime}v\right\rangle -\left\langle k,\wp_{-}^{\prime
}v\right\rangle \nonumber
\end{align}
defined properly in the subspace $\mathfrak{H}_{-}$ of the Hilbert space
(\ref{A12}). We call (\ref{A44}) an \textit{energy functional} for the problem
with point conditions.

\begin{remark}
\label{remENER}In the case $k=0$ the form $\mathfrak{E}(v;f,0)$ restricted
onto $\mathcal{D}(\mathcal{A})\times\mathcal{L}\subset\mathfrak{H}_{-}%
\times\mathcal{L}$ coincides with the energy functional%
\[
\frac{1}{2}(\mathcal{A}v,v)_{\mathcal{L}}-(f,v)_{\mathcal{L}}%
\]
generated by the self-adjoint extension $\mathcal{A}$ with the parameters
(\ref{A24}) in Proposition \ref{prop3A}. This follows from the fact that two
scalar products in $\mathbb{R}^{J}$ on the right-hand side of (\ref{A44}) vanish.
\end{remark}

\begin{theorem}
\label{th7A}A vector function $v\in\mathfrak{H}$ is a solution of the problem
(\ref{A27})-(\ref{A30}) if and only if $v$ is a stationary point of the energy
functional (\ref{A44}).
\end{theorem}

\textbf{Proof.} Calculating the variation of the functional (\ref{A44}), we
obtain%
\begin{align*}
\delta\mathfrak{E}(v,w;f,k)  &  =-\frac{1}{2}(\Delta_{y}v_{0},w_{0}%
)_{\omega_{0}}-\frac{1}{2}(\Delta_{y}w_{0},v_{0})_{\omega_{0}}-(f_{0}%
,v_{0})_{\omega_{0}}\\
&  -\sum\nolimits_{j}\left(  \frac{1}{2}\gamma_{j}|\omega_{j}|(\partial
_{z}^{2}v_{j},w_{j})_{I_{j}}+\frac{1}{2}\gamma_{j}|\omega_{j}|(\partial
_{z}^{2}w_{j},v_{j})_{I_{j}}+(f_{j},w_{j})_{I_{j}}\right) \\
&  +\frac{1}{2}\left\langle \wp_{+}^{\prime\prime}v-\wp_{+}^{\prime}v-S\wp
_{-}^{\prime}v,\wp_{-}^{\prime}w\right\rangle +\frac{1}{2}\left\langle \wp
_{+}^{\prime\prime}w-\wp_{+}^{\prime}w-S\wp_{-}^{\prime}w,\wp_{-}^{\prime
}v\right\rangle -\left\langle k,\wp_{-}^{\prime}w\right\rangle .
\end{align*}
We make use of the generalized Green formula (\ref{A14}) while interchanging
positions of $v$ and $w$. Recalling the relation $S=S^{\top}$ and the point
condition (\ref{A27}) for $v,w\in\mathfrak{H}_{-}$, we have%
\begin{align}
\delta\mathfrak{E}(v,w;f,k)  &  =-\frac{1}{2}(\Delta_{y}v_{0},w_{0}%
)_{\omega_{0}}-\frac{1}{2}(w_{0},\Delta_{y}v_{0})_{\omega_{0}}-(f_{0}%
,v_{0})_{\omega_{0}}\label{A45}\\
&  -\sum\nolimits_{j}\left(  \frac{1}{2}\gamma_{j}|\omega_{j}|(\partial
_{z}^{2}v_{j},w_{j})_{I_{j}}+\frac{1}{2}\gamma_{j}|\omega_{j}|(w_{j}%
,\partial_{z}^{2}v_{j})_{I_{j}}+(f_{j},w_{j})_{I_{j}}\right) \nonumber\\
&  +\frac{1}{2}\left\langle \wp_{+}w,\wp_{-}v\right\rangle +\frac{1}%
{2}\left\langle \wp_{-}w,\wp_{+}v\right\rangle +\frac{1}{2}\left\langle
\wp_{+}^{\prime\prime}v-\wp_{+}^{\prime}v,\wp_{-}^{\prime}w\right\rangle
\nonumber\\
&  +\frac{1}{2}\left\langle \wp_{+}^{\prime\prime}w-\wp_{+}^{\prime}w,\wp
_{-}^{\prime}v\right\rangle +\frac{1}{2}\left\langle S\wp_{-}^{\prime}%
v,\wp_{-}^{\prime}w\right\rangle -\frac{1}{2}\left\langle \wp_{-}^{\prime
}w,S\wp_{-}^{\prime}v\right\rangle -\left\langle k,\wp_{-}^{\prime
}w\right\rangle \nonumber\\
&  =(-\Delta_{y}v_{0}-f_{0},w_{0})_{\omega_{0}}+\sum\nolimits_{j}(-\gamma
_{j}|\omega_{j}|\partial_{z}^{2}v_{j}-f_{j},w_{j})_{I_{j}}+\left\langle
\wp_{+}^{\prime\prime}v-\wp_{+}^{\prime}v-S\wp_{-}^{\prime}v,\wp_{-}^{\prime
}w\right\rangle .\nonumber
\end{align}
Here we, in particular, used the relation (\ref{A27}). It also should be
mentioned that all functions are real as well as the matrix $S$ and the vector
$k$.

We see that a solution $v\in\mathfrak{H}$ of the problem (\ref{A27}%
)-(\ref{A30}), annuls the variation (\ref{A45}) of the functional (\ref{A44}).
On the other hand, for any test vector $w\in\mathfrak{H}_{-}$, the expression
with the stationary point $v\in\mathfrak{H}_{-}$ of $\mathfrak{E}$ vanishes;
in particular, taking $w\in C_{c}^{\infty}(\omega_{0})\times C_{c}^{\infty
}(I_{1})\times...\times C_{c}^{\infty}(I_{J})$ brings the differential
equations in (\ref{A29}) and (\ref{A30}). Thus, the last scalar product in
(\ref{A45}) is null and (\ref{A28}) is fulfilled because $\wp_{-}^{\prime
}\mathfrak{H}_{-}=\mathbb{R}^{J}.$ It remains to mention that the boundary
conditions (\ref{22}), (\ref{23}) and the point condition (\ref{A27}) are kept
in the space $\mathfrak{H}_{-}$. $\blacksquare$

\section{Determination of parameters of an appropriate hybrid
model\label{sect3}}

\subsection{The boundary layer phenomenon\label{sect3.1}}

An asymptotic analysis performed in \cite{BuCaNa1} gave a detailed description
of the behavior of solutions to the stationary problem in $\Xi(h)$ near the
junction zones. The internal constitution of the boundary layers which appear
in the vicinity of the sockets $\theta_{j}^{h}=\omega_{j}^{h}\times(0,h)$ and
are written in the rapid variables%
\begin{equation}
\xi^{j}=(\eta^{j},\zeta^{j}),\ \ \eta^{j}=h^{-1}(y-P^{j}),\ \ \zeta^{j}%
=h^{-1}z, \label{B1}%
\end{equation}
depends crucially on the exponent $\alpha$ in (\ref{14}). In the case
$\alpha=1$, see (\ref{15}), the transmission conditions (\ref{12}), (\ref{13})
decouple and the coordinate dilation leads to two independent limit problems
in the semi-infinite cylinder $Q_{j}=\omega_{j}\times\mathbb{R}^{+}$ and the
perforated layer $\Lambda_{j}=(\mathbb{R}^{2}\setminus\overline{\omega}%
_{j})\times(0,1)$. In \cite[Sect. 2.4]{BuCaNa1}, we have examined these
problems, namely the Neumann problem%
\begin{align}
-\gamma_{j}\Delta_{\xi}W_{j}\left(  \xi\right)   &  =0,\ \xi\in Q_{j}%
,\ \ \gamma_{j}\partial_{\nu}W_{j}\left(  \xi\right)  =g_{j}\left(
\xi\right)  ,\ \ \ \xi\in\partial\omega_{j}\times\mathbb{R}_{+},\label{B2}\\
-\gamma_{j}\partial_{\zeta}W_{j}\left(  \eta,0\right)   &  =0,\ \ \eta
\in\omega_{j},\nonumber
\end{align}
where $\partial_{\nu}$ is the outward normal derivative, and the mixed
boundary-value problem%
\begin{align}
-\Delta_{\xi}W_{0}\left(  \xi\right)   &  =0,\ \xi\in\Lambda_{j}%
,\ \ -\partial_{\zeta}W_{0}\left(  \eta,0\right)  =\partial_{\zeta}%
W_{0}\left(  \eta,1\right)  =0,\ \ \ \eta\in\mathbb{R}^{2}\setminus
\overline{\omega}_{j},\label{B3}\\
W_{0}\left(  \xi\right)   &  =g_{0}\left(  \xi\right)  ,\ \ \xi\in
\partial\omega_{j}\times(0,1).\nonumber
\end{align}

We now point out several special solutions of these problems, that we need in
the sequel. First of all, the homogeneous ($g_{j}=0$) problem (\ref{B2}) has a
constant solution, say $\mathbf{w}_{j}(\xi)=1$, and the problem (\ref{B3})
with $g_{0}(\xi)=1$ is also satisfied by $\mathbf{w}_{j}^{0}(\xi)=1$. The
homogeneous ($g_{0}=0$) problem (\ref{B3}) admits a solution with the
logarithmic growth at infinity%
\begin{equation}
\mathbf{W}_{j}^{0}\left(  \eta\right)  =(2\pi)^{-1}\left(  \ln|\eta|+\ln
c_{\log}\left(  \omega_{j}\right)  \right)  +\widetilde{\mathbf{W}}_{j}%
^{0}\left(  \eta\right)  ,\ \ \ \widetilde{\mathbf{W}}_{j}^{0}\left(
\eta\right)  =O\left(  |\eta|^{-1}\right)  ,\ \ |\xi|\rightarrow+\infty,
\label{B4}%
\end{equation}
where $c_{\log}\left(  \omega_{j}\right)  $ is the logarithmic capacity of the
set $\overline{\omega}_{j}\subset\mathbb{R}^{2}$. Note that the function
$\mathbf{W}_{j}^{0}$ in (\ref{B4}) is independent of $\zeta$ and is called the
logarithmic capacity potential, see \cite{PoSe, Land}. Finally, according to
\cite[Lemma 6]{BuCaNa1},%
\begin{equation}
-1=%
%TCIMACRO{\dint _{\partial\omega_{j}}}%
%BeginExpansion
{\displaystyle\int_{\partial\omega_{j}}}
%EndExpansion
g_{j}\left(  \eta\right)  ds_{\eta},\ \ \ g_{j}\left(  \eta\right)
=\partial_{\nu}\mathbf{W}_{j}^{0}\left(  \eta\right)  \label{B5}%
\end{equation}
and a solution of the Neumann problem (\ref{B2}) with the datum in (\ref{B5})
can be found in the form%
\begin{equation}
\mathbf{W}_{j}\left(  \eta,\zeta\right)  =\gamma_{j}^{-1}|\omega_{j}%
|^{-1}\zeta+O(e^{-\delta\zeta}),\ \ \delta>0. \label{B6}%
\end{equation}

\subsection{Individual choice of the self-adjoint extension\label{sect3.2}}

Applying the method of matched asymptotic expansions, see for example
\cite{VanDyke, Ilin}, \cite[Ch. 2]{MaNaPl}, we take some functions $v_{0}%
\in\mathfrak{H}_{0},$ $v_{j}\in\mathfrak{H}_{j}$ and write the \textit{outer}
expansions in the plate $\Omega_{0}(h)$ and the rod $\Omega_{j}(h)$ near but
outside the socket $\theta_{j}^{h}$%
\begin{align}
v_{0}(y)  &  =\frac{b_{j}}{2\pi}\ln\frac{1}{r_{j}}+\widehat{v}_{0}%
(P^{j})+...=\frac{b_{j}}{2\pi}\ln\frac{1}{|\eta_{j}|}-\frac{b_{j}}{2\pi}\ln
h+\widehat{v}_{0}(P^{j})+...,\label{B7}\\
v_{j}(z)  &  =v_{j}(0)+z\partial_{z}v_{j}(0)+...=v_{j}(0)+h\zeta^{j}%
v_{j}(0)+... \label{B8}%
\end{align}
Here, ellipses stand for higher order terms of no importance in our asymptotic
procedure. The \textit{inner} expansions in the immediate vicinity of the
sockets are composed from the above described solutions of the problems
(\ref{B3}) and (\ref{B2})%
\begin{align}
c_{j}^{0}\mathbf{w}_{j}^{0}(\xi^{j})+c_{j}^{1}\mathbf{W}_{j}^{0}(\xi^{j})  &
=c_{j}^{1}\frac{1}{2\pi}\ln\frac{1}{|\eta_{j}|}+c_{j}^{1}\frac{1}{2\pi}\ln
c_{\log}\left(  \omega_{j}\right)  +c_{j}^{0}+...,\ (P^{j}+h\eta^{j}%
,h\zeta^{j})\in\Omega_{\bullet}(h),\label{B9}\\
c_{j}^{0}\mathbf{w}_{j}(\xi^{j})+hc_{j}^{1}\mathbf{W}_{j}(\xi^{j})  &
=hc_{j}^{1}\gamma_{j}^{-1}|\omega_{j}|^{-1}\zeta^{j}+c_{j}^{0}%
+...,\ \ \ (P^{j}+h\eta^{j},h\zeta^{j})\in\Omega_{j}(h). \label{B10}%
\end{align}
We emphasize that, by definition of those solutions, the transmission
condition (\ref{13}) is wholly satisfied by (\ref{B9}) and (\ref{B10}) but the
condition (\ref{13}) with the reasonable precision $O(h).$

Comparing (\ref{B7}) with (\ref{B9}) and (\ref{B8}) with (\ref{B10}) yields
the equations%
\begin{align}
b_{j}  &  =c_{j}^{1},\ \ \widehat{v}_{0}(P^{j})-b_{j}(2\pi)^{-1}\ln
h=c_{j}^{0}+c_{j}^{1}(2\pi)^{-1}\ln c_{\log}\left(  \omega_{j}\right)
,\label{BBB}\\
\partial_{z}v_{j}(0)  &  =c_{j}^{1}\gamma_{j}^{-1}|\omega_{j}|^{-1}%
,\ \ v_{j}(0)=c_{j}^{0}.\nonumber
\end{align}
Excluding the coefficients $c_{j}^{1}$ and $c_{j}^{0}$, we derive the
relations%
\begin{align}
b_{j}-\gamma_{j}|\omega_{j}|\partial_{z}v_{j}(0)  &  =0,\label{B11}\\
\widehat{v}_{0}(P^{j})-v_{j}(0)  &  =b_{j}(2\pi)^{-1}(\ln h+\ln c_{\log
}\left(  \omega_{j}\right)  ).\nonumber
\end{align}
We also introduce the diagonal $J\times J$-matrix%
\begin{equation}
S^{h}=-(2\pi)^{-1}\mathrm{diag}\{\ln h+\ln c_{\log}\left(  \omega_{1}\right)
,...,\ln h+\ln c_{\log}\left(  \omega_{J}\right)  \}=(2\pi)^{-1}(|\ln
h|\mathbb{I}+C_{\log}) \label{B12}%
\end{equation}
which is positive definite for $h\in(0,h_{0})$ with a small $h_{0}\in(0,1)$
and further substitutes for $S$ in (\ref{A24}). Here, $\mathbb{I}$ is the unit
$J\times J$-matrix and $C_{\log}=\mathrm{diag}\{\ln c_{\log}\left(  \omega
_{1}\right)  ,...,\ln c_{\log}\left(  \omega_{J}\right)  \}$.

The main conclusion from the above consideration is that the relations
(\ref{B11}) between the functions $v_{0}\in\mathfrak{H}_{0}$ and $v_{j}%
\in\mathfrak{H}_{j}$, $J=1,...,J$ are equivalent to conditions in (\ref{A20}),
(\ref{A24}) defining a self-adjoint extension $\mathcal{A}^{h}$ of the
operator $A=(A_{0},A_{1},...,A_{J})$ in the space $\mathfrak{L}$, see
(\ref{LL}), examined in Section \ref{sect2}. In what follows we deal with this
operator $\mathcal{A}^{h}$.

\subsection{The spectral problem\label{sect3.3}}

Based on results in Section \ref{sect2}, we give two models of the spectral
problem (\ref{7})-(\ref{13}). First, we use the introduced self-adjoint
operator $\mathcal{A}^{h}$ and write the equation%
\begin{equation}
\mathcal{A}^{h}v^{h}=\mu^{h}v^{h}\text{ \ in }\mathfrak{L}. \label{S1}%
\end{equation}
Second, we supply the problems (\ref{20}), (\ref{22}) and (\ref{21}),
(\ref{23}) with the point conditions (\ref{A26}), (\ref{A27}) and formulate
the abstract equation%
\begin{equation}
\mathfrak{A}^{h}v^{h}=\mu^{h}(v^{h},0)\text{ \ in }\mathfrak{L\times
}\mathbb{R}^{J} \label{S2}%
\end{equation}
where the operator $\mathfrak{A}^{h}:\mathfrak{H}_{-}\rightarrow
\mathfrak{L\times}\mathbb{R}^{J}$ involves the point condition (\ref{A26})
with the operator $S=S^{h}$ in (\ref{B12}) while zero at the last position in
(\ref{S2}) indicates that this condition is homogeneous, namely $k=0$ in
(\ref{A28}).

According to Remark \ref{remEQU}, the spectral equations (\ref{S1}) and
(\ref{S2}) are equivalent with each other. Unfortunately, the operator
$\mathcal{A}^{h}$ intended to model the problem (\ref{7})-(\ref{13}) with the
positive spectrum (\ref{18}), is not positive definite in view of Remark
\ref{remENER} and formula (\ref{A44}) involving the negative definite matrix
$-S=-S^{h}$. However, in Remark \ref{remSIMPLE} we mentioned the positive
operator $\mathcal{A}^{0}$ with the attributes (\ref{A21}) in (\ref{A20})
which differs from the operator $\mathcal{A}^{h}$ only in subspace of
dimension $J$. Thus, the max-min principle, cf. \cite[Thm. 10.2.2]{BiSo}
demonstrates that the total multiplicity of the negative part of the spectrum
of $\mathcal{A}^{h}$ cannot exceed $J$. In the next section we will construct
asymptotics of negative eigenvalues which we call \textit{parasitic}.

\subsection{Asymptotics of parasitic eigenvalues\label{sect3.4}}

We will need the fundamental solution $\Phi$ of the operator $-\Delta_{y}+1$
in the plane $\mathbb{R}^{2}$. Its expression is well known, in particular%
\begin{equation}
\Phi(y)=O(e^{-\psi|y|})\text{ as }|y|\rightarrow+\infty,\ \ \ \Phi(y)=\frac
{1}{2\pi}\ln\frac{1}{|y|}+\Psi+O(|y|)\text{ as }|y|\rightarrow+0 \label{S3}%
\end{equation}
but exact values of $\psi>0$ and $\Psi$ are of no importance here. Let%
\begin{equation}
\mu^{h}=-h^{-2}(e^{2m}+O(h)). \label{S4}%
\end{equation}
We set%
\begin{align}
v_{0}^{h}(y)  &  =%
%TCIMACRO{\tsum \nolimits_{j}}%
%BeginExpansion
{\textstyle\sum\nolimits_{j}}
%EndExpansion
\alpha_{j}\chi_{j}(y)\Phi(h^{-1}e^{m}(y-P^{j})),\label{S5}\\
v_{j}^{h}(y)  &  =\alpha_{j}X_{j}(z)\gamma_{j}^{-1}|\omega_{j}|^{-1}%
he^{-m}e^{-h^{-1}e^{m}z}. \label{S6}%
\end{align}
where the column $\alpha=(\alpha_{1},...,\alpha_{J})\in\mathbb{R}^{J}$ and the
number $m\in\mathbb{R}$ are to be determined and $X_{j}\in C^{\infty}(I_{j}),$
$X_{j}(z)=1$ for $[0,l_{j}/3]$ and $X_{j}(z)=0$ for $[2l_{j}/3,l_{j}]$.
Clearly, the functions (\ref{S5}), (\ref{S6}) satisfy the boundary conditions
(\ref{22}), (\ref{23}) and leave small discrepancies $O(e^{-\delta/h}),$
$\delta>0,$ in the differential equations (\ref{20}), (\ref{21}) with the
spectral parameter $-h^{-2}e^{2m}$. It should be mentioned that the
exponential decay of the boundary layer terms in (\ref{S5}) and (\ref{S6}) is
due to the negative value of $\mu^{h}$ in (\ref{S4}).

The vector $v^{h}=(v_{0}^{h},v_{1}^{h},...,v_{J}^{h})$ has the following
projections:%
\begin{align}
\wp_{+}^{\prime}v^{h}  &  =(2\pi)^{-1}(m-\ln h)+\phi_{0},\ \ \wp_{-}^{\prime
}v^{h}=\alpha,\label{S7}\\
\wp_{+}^{\prime\prime}v^{h}  &  =Qhe^{-m},\ \ \wp_{-}^{\prime\prime}%
v^{h}=-\alpha,\nonumber
\end{align}
where $Q=\mathrm{diag}\{\gamma_{1}^{-1}|\omega_{1}|^{-1},...,\gamma_{J}%
^{-1}|\omega_{J}|^{-1}\}.$ Hence, the point condition (\ref{A27}) is fulfilled
while, in view of (\ref{B12}) and (\ref{S7}), the condition (\ref{A26})
converts into%
\[
Qhe^{m}\alpha-(2\pi)^{-1}(m-\ln h)\alpha+(2\pi)^{-1}(\ln h\mathbb{I}+C_{\log
})\alpha=0
\]
or, what is the same,%
\begin{equation}
C_{\log}\alpha=m\alpha-2\pi he^{-m}Q\alpha. \label{S8}%
\end{equation}
Since both matrices are diagonal, the system (\ref{S8}) splits into $J$
independent transcendental equations. The small factor $h$ of the exponent
$e^{-m}$ allows us to apply the implicit function theorem and to obtain the
solutions%
\begin{equation}
m_{j}^{h}=\log c_{\log}(\omega_{j})+O(h),\ \ \ \alpha_{(j)}^{h}=\mathbf{e}%
_{(j)}+O(h),\ \ \ j=1,...,J. \label{S9}%
\end{equation}

We have derived "good approximations" for $J$ negative eigenvalues of the
spectral equations (\ref{S1}) and (\ref{S2}). Recalling that their number
cannot exceed $J,$ we come in position to formulate an assertion on the whole
negative part of the spectrum.

\begin{proposition}
\label{prop11S}There exist positive numbers $h_{-}$ and $c_{-}$ such that, for
$h\in(0,h_{-}]$, the equation (\ref{S1}) or (\ref{S2}) possesses exactly $J$
eigenvalues on the semi-axis $\mathbb{R}_{-}=(-\infty,0).$ These eigenvalues
obey the asymptotic form%
\begin{equation}
|\mu_{-j}^{h}+h^{-2}(c_{\log}(\omega_{j}))^{2}|\leq c_{-}h^{-1} \label{S10}%
\end{equation}
where $c_{\log}(\omega_{j})>0$ is the logarithmic capacity of the set
$\overline{\omega}_{j}\subset\mathbb{R}^{2}$, see \cite{PoSe, Land}.
\end{proposition}

We will outline the proof in Remark \ref{remNEGA}. Notice that the negative
eigenvalues $\mu_{-j}^{h}=O(h^{-2}),$ $j=1,...,J,$ are situated very far away
from the positive part of the spectrum, that is, outside the scope of the
asymptotic models under consideration.

\section{Justification of the asymptotic models\label{sect4}}

\subsection{The first convergence theorem\label{sect4.1}}

Let $\mu^{h}$ be an eigenvalue of the self-adjoint operator $\mathcal{A}^{h}$
in (\ref{S1}). The corresponding eigenvector $v^{h}=(v_{0}^{h},v_{1}%
^{h},...,v_{J}^{h})\in\mathcal{D}(\mathcal{A}^{h})\subset\mathfrak{H}_{-}$ can
be normed as follows:%
\begin{equation}
||v_{0}^{h};L^{2}(\omega_{0})||^{2}+%
%TCIMACRO{\tsum \nolimits_{j}}%
%BeginExpansion
{\textstyle\sum\nolimits_{j}}
%EndExpansion
\rho_{j}|\omega_{j}|~||v_{j}^{h};L^{2}(I_{j})||^{2} \label{T1}%
\end{equation}
Assuming that%
\begin{equation}
|\mu^{h}|\leq c, \label{T2}%
\end{equation}
in particular, rejecting negative eigenvalues in Proposition \ref{prop11S}, we
recall the Kondratiev theory used in Section \ref{sect2.1}. Then, a solution
$v_{0}^{h}\in H_{loc}^{2}(\overline{\omega}_{0}\setminus\mathcal{P})\cap
L^{2}(\omega_{0})$ of the problem (\ref{20}), (\ref{22}) admits the
decomposition (\ref{A9}) with the ingredients $a_{j}^{h}=\widehat{v}^{h}%
(P^{j}),$ $b_{j}^{h}\in\mathbb{R}$ and $\widehat{v}_{0}^{h}\in H^{2}%
(\omega_{0})\cap H_{0}^{1}(\omega_{0})$ while%
\begin{equation}
||\widehat{v}_{0}^{h};H^{2}(\omega_{0})||^{2}+%
%TCIMACRO{\tsum \nolimits_{j}}%
%BeginExpansion
{\textstyle\sum\nolimits_{j}}
%EndExpansion
|b_{j}^{h}|\leq c(1+|\mu^{h}|)||v_{0}^{h};L^{2}(\omega_{0})||^{2}\leq C,
\label{T3}%
\end{equation}
where $\widehat{v}_{0}^{h}$ is given in (\ref{A99}).

Furthermore, a solution $v_{j}^{h}\in L^{2}(I_{j})$ of the ordinary
differential equation (\ref{21}) with the Dirichlet condition (\ref{23}) falls
into the Sobolev space $H^{2}(I_{j})$ and fulfills the estimate%
\begin{equation}
|v_{j}^{h}(0)|+|\partial_{z}v_{j}^{h}(0)|\leq c||v_{j}^{h};H^{2}(I_{j})||\leq
c(1+|\mu^{h}|)||v_{0}^{h};L^{2}(I_{j})||^{2}\leq C. \label{T4}%
\end{equation}

The inequalities (\ref{T2})-(\ref{T4}) help to conclude the following
convergence along an infinitesimal positive sequence $\{h_{k}\}_{k\in
\mathbb{N}}$:%
\begin{align}
\mu^{h}  &  \rightarrow\mu^{0}\in\mathbb{R},\ \ \ b^{h}=(b_{1}^{h}%
,...,b_{J}^{h})\rightarrow b^{0}\in\mathbb{R}^{J},\ \label{T5}\\
\widehat{v}_{0}^{h}  &  \rightharpoonup\widehat{v}_{0}^{0}\text{ weakly in
}H^{2}(\omega_{0}),\ v_{j}^{h}\rightharpoonup v_{j}^{0}\text{ weakly in }%
H^{2}(I_{j}).\nonumber
\end{align}
This and the embeddings $H^{2}(\omega_{0})\subset C(\omega_{0}),$ $H^{2}%
(I_{j})\subset C^{1}(I_{j})$ imply the convergence of the projections
(\ref{A13})%
\begin{equation}
\wp_{\pm}v^{h}\rightarrow\wp_{\pm}v^{0}\in\mathbb{R}^{J}. \label{T6}%
\end{equation}
We emphasize that formulas in (\ref{T5}) guarantee the strong convergences
$v_{0}^{h}\rightarrow v_{0}^{0}$ in $L^{2}(\omega_{0})$ and $v_{j}%
^{h}\rightarrow v_{j}^{0}$ in $L^{2}(I_{j})$ so that the normalization
condition (\ref{T1}) is kept by the limit $v^{0}=(v_{0}^{0},v_{1}%
^{0},...,v_{J}^{0})$. Moreover, the differential equations (\ref{20}),
(\ref{21}) with $\mu=\mu^{h}$ and the boundary conditions (\ref{22}),
(\ref{23}) for $v_{0}^{h},$ $v_{j}^{h}$ are passed to the limits%
\begin{equation}
\mu^{0},\ \ \ v_{0}^{0}=\widehat{v}_{0}^{0}-(2\pi)^{-1}%
%TCIMACRO{\tsum \nolimits_{j}}%
%BeginExpansion
{\textstyle\sum\nolimits_{j}}
%EndExpansion
\chi_{j}b_{j}\ln r_{j},\text{ \ \ }v_{j}^{h},\ j=1,...,J. \label{T7}%
\end{equation}

In order to formulate the next assertion it suffices to mention that the point
conditions (\ref{A26}), (\ref{A27}) with the matrix (\ref{B12}) containing the
big component $-(2\pi)^{-1}\ln h\mathbb{I}$ turn in the limit into the
relations%
\begin{equation}
\wp_{-}^{\prime}v^{0}=0,\ \ \wp_{-}^{\prime}v^{0}+\wp_{-}^{\prime\prime}%
v^{0}=0\ \ \Rightarrow\ \ \wp_{-}^{\prime}v^{0}=\wp_{-}^{\prime\prime}%
v^{0}=0\in\mathbb{R}^{J}. \label{T8}%
\end{equation}
These provide the self-adjoint extension $\mathcal{A}^{0}$ with the attributes
(\ref{A21}) in (\ref{A20}) that corresponds to the Neumann and mixed
boundary-value problems (\ref{A22}) and (\ref{A23}). The spectra
$\{\varkappa_{n}^{0}\}_{n\in\mathbb{N}}$ and $\{\varkappa_{n}^{j}=\frac
{\pi^{2}}{l_{j}^{2}}\frac{\gamma_{j}}{\rho_{j}}(n+\frac{1}{2})^{2}%
\}_{n\in\mathbb{N}},$ $j=1,...,J,$ of the above mentioned problems are united
into the common monotone sequence%
\begin{equation}
\{\mu_{n}^{0}\}_{n\in\mathbb{N}},\ \ \mu_{1}^{0}=0<\mu_{2.}^{0} \label{T9}%
\end{equation}

Here, eigenvalues are listed while counting their multiplicity in.

\begin{theorem}
\label{th1T}If an eigenvalue $\mu^{h}$ of the operator $\mathcal{A}^{h}$, cf.
(\ref{S1}) and (\ref{S2}), and the corresponding eigenvector $v^{h}$ fulfil
the requirements (\ref{T2}) and (\ref{T1}), then the limits $\mu^{0}$ and
$v^{0}$ in (\ref{T5}) along an infinitesimal sequence $\{h_{n}\}_{n\in
\mathbb{N}}$ are an eigenvalue of the operator $\mathcal{A}^{0}$ described in
Remark \ref{remSIMPLE} and the corresponding eigenvector normed in the space
$\mathcal{L}$, see (\ref{LL}).
\end{theorem}

\subsection{The second convergence theorem\label{sect4.2}}

In the next section we will verify that entries of the eigenvalue sequence
(\ref{18}) of the original problem (\ref{7})-(\ref{13}) in the junction
$\Xi(h)\subset\mathbb{R}^{3}$ satisfy the inequalities%
\begin{equation}
0<\lambda^{n}(h)\leq c_{n}\text{ for }h\in(0,h_{n}] \label{R1}%
\end{equation}
with some positive $h_{n}$ and $c_{n}$ which depend on the eigenvalue number
$n$ but are independent of $h$. The corresponding eigenfunction $u^{n}%
(h,\cdot)\in H_{0}^{1}(\Xi(h);\Gamma(h))$ is subject to the normalization
condition (\ref{19}). We introduce the functions%
\begin{align}
v_{0}^{n}(h,y)  &  =\frac{1}{\sqrt{h}}\int_{0}^{h}u^{n}(h,y,z)dz,\ \ y\in
\omega_{0},\label{R2}\\
v_{j}^{n}(h,z)  &  =\frac{1}{h^{3/2}|\omega_{j}|}\int_{\omega_{j}^{h}}%
u_{j}^{n}(h,y,z)dy,\ \ y\in(0,l_{j}),\ j=1,...,J, \label{R3}%
\end{align}
and write%
\begin{align}
\int_{\omega_{0}}|v_{0}^{n}(h,y)|^{2}dy  &  =\frac{1}{h}\int_{\omega_{0}%
}\left\vert \int_{0}^{h}u^{n}(h,y,z)dz\right\vert ^{2}dy\leq\int_{\omega_{0}%
}\int_{0}^{h}\left\vert u^{n}(h,x)\right\vert ^{2}dx\label{R4}\\
&  \leq b(u^{n},u^{n};\Omega_{0}(h)),\nonumber\\
\rho_{j}|\omega_{j}|\int_{0}^{l_{j}}|v_{j}^{n}(h,z)|^{2}dz  &  =\rho
_{j}|\omega_{j}|^{-1}\frac{1}{h^{3}}\int_{0}^{l_{j}}\left\vert \int
_{\omega_{j}^{h}}u_{j}^{n}(h,y,z)dy\right\vert ^{2}dz\nonumber\\
&  \leq\frac{\rho_{j}}{h}\frac{|\omega_{j}^{h}|}{h^{2}|\omega_{j}|}%
\int_{\Omega_{j}^{h}(h)}|u^{n}(h,x)|^{2}dx\leq b(u^{n},u^{n};\Omega
_{j}(h)).\nonumber
\end{align}
Here, we used formulas (\ref{17}) and (\ref{14}), (\ref{15}) while taking the
relation $1\leq h^{-1}\rho_{j}$ on $\Omega_{0}(h)\cap\Omega_{j}(h)$ into
account. Hence, the vector function $v^{n}=(v_{0}^{n},v_{1}^{n},...,v_{J}%
^{n})$ satisfies the estimate%
\begin{equation}
||v^{n}(h,\cdot);\mathfrak{L}||\leq1. \label{R5}%
\end{equation}
A similar calculation gives us the formula%
\begin{align*}
||\nabla v_{0}^{n};L^{2}(\omega_{0})||^{2}+%
%TCIMACRO{\tsum \nolimits_{j}}%
%BeginExpansion
{\textstyle\sum\nolimits_{j}}
%EndExpansion
||\partial_{z}v_{j}^{n};L^{2}(I_{j})||^{2}  &  \leq c(||\nabla_{y}u^{n}%
;L^{2}(\Omega_{0}(h))||^{2}+%
%TCIMACRO{\tsum \nolimits_{j}}%
%BeginExpansion
{\textstyle\sum\nolimits_{j}}
%EndExpansion
||\partial_{z}u^{n};L^{2}(\Omega_{j}(h))||^{2})\\
&  \leq ca(u^{n},u^{n};\Xi(h))=c\lambda(h)b(u^{n},u^{n};\Xi(h))\leq C_{n}.
\end{align*}

Moreover, the Poincar\'{e} inequalities in $(0,h)$ and $\omega_{j}^{h}$ show
that the functions%
\begin{align*}
u_{0}^{n\perp}(h,x)  &  =u^{n}(h,x)-h^{-1/2}v_{0}^{n}(h,y)\text{\ \ in }%
\Omega_{0}(h),\\
u_{j}^{n\perp}(h,x)  &  =u_{j}^{n}(h,x)-h^{-1/2}v_{j}^{n}(h,y)\text{\ \ in
}\Omega_{j}(h),
\end{align*}
which are of mean zero in $z\in(0,h)$ and $y\in\omega_{j}^{h}$, respectively,
enjoy the relations%
\begin{align}
||u_{0}^{n\perp};L^{2}(\Omega_{0}(h))||^{2}  &  \leq c_{0}h^{2}||\partial
_{z}(u^{n}-h^{-1/2}v_{0}^{n});L^{2}(\Omega_{0}(h))||^{2}\label{R6}\\
&  =c_{0}h^{2}||\partial_{z}u^{n};L^{2}(\Omega_{0}(h))||^{2}\leq c_{0}%
a(u^{n},u^{n};\Omega_{0}(h)),\nonumber\\
||u_{j}^{n\perp};L^{2}(\Omega_{j}(h))||^{2}  &  \leq c_{j}h^{2}||\nabla
_{y}(u_{j}^{n}-h^{-1/2}v_{j}^{n});L^{2}(\Omega_{j}(h))||^{2}\nonumber\\
&  =c_{j}h^{2}||\nabla_{y}u^{n};L^{2}(\Omega_{j}(h))||^{2}\leq C_{j}%
h^{3}a(u^{n},u^{n};\Omega_{j}(h)).\nonumber
\end{align}
Hence, we obtain%
\begin{align}
1  &  =||h^{-1/2}v_{0}^{n}-u_{0}^{n\perp};L^{2}(\Omega_{\bullet}%
(h))||^{2}+h^{-1}%
%TCIMACRO{\tsum \nolimits_{j}}%
%BeginExpansion
{\textstyle\sum\nolimits_{j}}
%EndExpansion
\rho_{j}||h^{-1/2}v_{j}^{n}-u_{j}^{n\perp};L^{2}(\Omega_{j}(h))||^{2}%
\label{R7}\\
&  \leq(1+h)(||v_{0}^{n};L^{2}(\omega_{\bullet}(h))||^{2}+%
%TCIMACRO{\tsum \nolimits_{j}}%
%BeginExpansion
{\textstyle\sum\nolimits_{j}}
%EndExpansion
\rho_{j}|\omega_{j}|~||v_{j}^{n};L^{2}(I_{j})||^{2})\nonumber\\
&  +(1+h^{-1})(||u_{0}^{n\perp};L^{2}(\Omega_{\bullet}(h))||^{2}+h^{-1}%
%TCIMACRO{\tsum \nolimits_{j}}%
%BeginExpansion
{\textstyle\sum\nolimits_{j}}
%EndExpansion
\rho_{j}||u_{j}^{n\perp};L^{2}(\Omega_{j}(h))||^{2}).\nonumber
\end{align}
To estimate the norm $||v_{0}^{n};L^{2}(\omega_{\bullet}(h))||^{2}$, we apply
the weighted inequality in \cite[Thm. 9]{BuCaNa1}%
\begin{align}
&  (1+|\ln h|)^{-2}||r^{-1}(1+|\ln r|)^{-1}u;L^{2}(\Omega_{0}(h))||^{2}+h^{-1}%
%TCIMACRO{\tsum \nolimits_{j}}%
%BeginExpansion
{\textstyle\sum\nolimits_{j}}
%EndExpansion
||(l_{j}-z)^{-1}u_{j};L^{2}(\Omega_{j}(h))||^{2}\label{R8}\\
&  \leq c_{\Xi}a(u,u;\Xi(h))\nonumber
\end{align}
where $r=\min\{r_{1},...,r_{J}\}$ and $c_{\Xi}$ is independent of
$h\in(0,h_{0}]$ and $u\in H_{0}^{1}(\Xi(h);\Gamma(h))$. Since $r_{j}(1+|\ln
r_{j}|)\leq c_{j}h(1+|\ln h|)$ in the socket $\theta_{j}^{h}=\omega_{j}%
^{h}\times(0,h)$, we recall (\ref{16}), (\ref{R1}) and conclude that%
\begin{align}
||v_{0}  &  ^{n};L^{2}(\omega_{j}^{h})||^{2}\leq||u^{n};L^{2}(\theta_{j}%
^{h})||^{2}\leq ch^{2}(1+|\ln h|)^{4}||r_{j}^{-1}(1+|\ln r_{j}|)^{-1}%
u^{n};L^{2}(\theta_{j}^{h})||^{2}\label{R9}\\
&  \leq ch^{2}(1+|\ln h|)^{4}a(u^{n},u^{n};\Xi(h))=ch^{2}(1+|\ln
h|)^{4}\lambda^{n}(h)b(u^{n},u^{n};\Xi(h))\leq C_{n}h^{2}(1+|\ln
h|)^{4}.\nonumber
\end{align}
The estimates (\ref{R1}) and (\ref{R4}), (\ref{R5}) provide the following
convergence along an infinitesimal sequence $\{h_{k}\}_{k\in\mathbb{N}}$:%
\begin{align}
\lambda^{n}(h)  &  \rightarrow\lambda^{0},\label{R10}\\
v_{0}^{n}(h,\cdot)  &  \rightharpoondown v_{0}^{n0}\text{ weakly in }%
H^{1}(\omega_{0})\text{ and strongly in }L^{2}(\omega_{0}),\nonumber\\
v_{j}^{n}(h,\cdot)  &  \rightharpoondown v_{j}^{n0}\text{ weakly in }%
H^{1}(I_{j})\text{ and strongly in }L^{2}(I_{j}).\nonumber
\end{align}
We compose a test function $w$ from components $w_{0}\in C_{c}^{\infty
}(\overline{\omega}_{0}\setminus\mathcal{P})$ and $w_{j}\in C_{c}^{\infty
}(I_{j}),$ $j=1,...,J$. Then according to (\ref{17}) and (\ref{R2}),
(\ref{R3}), we transform the integral identity into the formula%
\begin{align*}
\sqrt{h}(\nabla_{y}  &  v_{0}^{n},\nabla_{y}w_{0})_{\omega_{0}}+\sqrt{h}%
%TCIMACRO{\tsum \nolimits_{j}}%
%BeginExpansion
{\textstyle\sum\nolimits_{j}}
%EndExpansion
\gamma_{j}|\omega_{j}|(\partial_{z}v_{j}^{n},\partial_{z}w_{j})_{I_{j}%
}=(\nabla_{y}u^{n},\nabla_{y}w_{0})_{\Omega_{\bullet}(h)}+h^{-1}%
%TCIMACRO{\tsum \nolimits_{j}}%
%BeginExpansion
{\textstyle\sum\nolimits_{j}}
%EndExpansion
\gamma_{j}(\partial_{z}u_{j}^{n},\partial_{z}w_{j})_{\Omega_{j}(h)}\\
&  =\lambda(h)\left(  (u^{n},w_{0})_{\Omega_{\bullet}(h)}+h^{-1}%
%TCIMACRO{\tsum \nolimits_{j}}%
%BeginExpansion
{\textstyle\sum\nolimits_{j}}
%EndExpansion
\rho_{j}(u_{j}^{n},w_{j})_{\Omega_{j}(h)}\right)  =\sqrt{h}\lambda(h)\left(
(v_{0}^{n},w_{0})_{\omega_{0}}+%
%TCIMACRO{\tsum \nolimits_{j}}%
%BeginExpansion
{\textstyle\sum\nolimits_{j}}
%EndExpansion
\rho_{j}|\omega_{j}|(v_{j}^{n},w_{j})_{I_{j}}\right)  .
\end{align*}
We multiply this with $h^{-1/2}$ and perform the limit passage along the
sequence $\{h_{k}\}_{k\in\mathbb{N}}$ to obtain%
\begin{equation}
(\nabla_{y}v_{0}^{n0},\nabla_{y}w_{0})_{\omega_{0}}+%
%TCIMACRO{\tsum \nolimits_{j}}%
%BeginExpansion
{\textstyle\sum\nolimits_{j}}
%EndExpansion
\gamma_{j}|\omega_{j}|(\partial_{z}v_{j}^{n0},\partial_{z}w_{j})_{I_{j}%
}=\lambda^{0}\left(  (v_{0}^{n0},w_{0})_{\omega_{0}}+%
%TCIMACRO{\tsum \nolimits_{j}}%
%BeginExpansion
{\textstyle\sum\nolimits_{j}}
%EndExpansion
\rho_{j}|\omega_{j}|(v_{j}^{n0},w_{j})_{I_{j}}\right)  . \label{R11}%
\end{equation}
Thanks to \cite[\S 9 Ch. 2]{LiMa}, the weak solution $v_{0}^{n0}\in
H^{1}(\omega_{0})$ of the variational problem (\ref{R11}) where $w_{j}=0,$
$j=1,...,J$, falls into $H^{2}(\omega_{0})$ and satisfies the Neumann problem
(\ref{20}), (\ref{22}) with $\mu=\lambda^{0}$. We emphasize that
$C_{c}^{\infty}(\overline{\omega}_{0}\setminus\mathcal{P})$ is dense in
$H^{1}(\omega_{0})$ so that any test function $w\in H^{1}(\omega_{0})$ is
available. At the same time, $w_{j}\in C_{c}^{\infty}(I_{j})$ vanishes near
the points $z=l_{j}$ and $z=0$. Hence, we may conclude that $v_{j}^{n0}\in
H^{2}(I_{j})$ and the differential equation (\ref{21}) with $\mu=\lambda^{0}$.
However, the boundary conditions%
\begin{equation}
v_{j}^{n0}(l_{j})=0,\ -\gamma_{j}|\omega_{j}|\partial_{z}v_{j}^{n0}(0)=0
\label{R12}%
\end{equation}
still must be derived. The Dirichlet condition in (\ref{R12}) is inherited
from the conditions $v_{j}^{n}(h,l_{j})=0$ and $u_{j}^{n}(h,y,l_{j})=0,$
$y\in\omega_{j}^{h}$. To conclude with the Neumann condition, we observe that
the inequality (\ref{R9}) allows us to repeat the above transformations with
the "very special" test vector function%
\begin{equation}
w_{0}^{j}(y)=\chi_{j}(y),\ \ w_{k}^{j}(z)=X_{k}(z)\delta_{k,j}%
,\ \ j,k=1,...,J, \label{R13}%
\end{equation}
where $\chi_{j}$ and $X_{k}$ are taken from (\ref{A5}) and (\ref{S6}). As a
result, the obtained information on $v_{0}^{n0}$ and $v_{j}^{n0}$ reduces the
integral identity (\ref{R12}) with (\ref{R13}) to the formula%
\[
0=-((\Delta_{y}+\lambda^{0})v_{0}^{n0},w_{0}^{j})_{\omega_{0}}-\gamma
_{j}|\omega_{j}|(((\partial_{z}^{2}+\lambda^{0})v_{j}^{n0},\partial_{z}%
w_{j}^{j})_{I_{j}}+\partial_{z}v_{j}^{n0}(0)w_{j}^{j}(0))=-\gamma_{j}%
|\omega_{j}|\partial_{z}v_{j}^{n0}(0).
\]
We are in position to formulate the convergence theorem.

\begin{theorem}
\label{th2T}The limits $\lambda^{0}$ and $v^{n0}=(v_{0}^{n0},v_{1}%
^{n0},...,v_{J}^{n0})$ in (\ref{R10}) are an eigenvalue and the corresponding
eigenvector normed by (\ref{T1}) of the problems (\ref{20}), (\ref{22}) and
(\ref{21}), (\ref{R12}).
\end{theorem}

\subsection{An abstract formulation of the original problem\label{sect4.3}}

In the Hilbert space $H^{h}=H_{0}^{1}(\Xi(h);\Gamma(h))$ we introduce the
scalar product%
\begin{equation}
(u^{h},v^{h})_{h}=a(u^{h},v^{h};\Xi(h))+b(u^{h},v^{h};\Xi(h)) \label{E1}%
\end{equation}
and positive, symmetric and continuous, therefore, self-adjoint operator
$T^{h}$,%
\begin{equation}
(T^{h}u^{h},v^{h})_{h}=b(u^{h},v^{h};\Xi(h))\ \ \ \ \forall u^{h},v^{h}\in
H^{h}. \label{E2}%
\end{equation}
Bilinear form on the right-hand side of (\ref{E1}) are defined in (\ref{17}).
Comparing (\ref{16}) with (\ref{E1}), (\ref{E2}), we see that the variational
formulation of the problem (\ref{7})-(\ref{13}) in $\Xi(h)$ is equivalent to
the abstract equation%
\begin{equation}
T^{h}u^{h}=\tau(h)u^{h}\text{ \ in }H^{h} \label{E3}%
\end{equation}
with the new spectral parameter%
\begin{equation}
\tau(h)=(1+\lambda(h))^{-1}. \label{E4}%
\end{equation}
The operator $T^{h}$ is compact and, hence, the essential spectrum of $T^{h}$
consists of the only point $\tau=0$, see, e.g., \cite[Thm. 10.1.5]{BiSo},
while the discrete spectrum composes the positive monotone infinitesimal
sequence%
\begin{equation}
1>\tau^{1}(h)>\tau^{2}(h)\geq...\geq\tau^{n}(h)\geq....\rightarrow+0
\label{E5}%
\end{equation}
obtained from (\ref{18}) according to formula (\ref{E4}).

The following assertion is known as the lemma on "near eigenvalues and
eigenvectors" \cite{ViLu} following directly from the spectral decomposition
of resolvent, see, e.g., \cite[Ch. 6]{BiSo}.

\begin{lemma}
\label{lemVL}Let $\mathcal{U}^{h}\in H^{h}$ and $t^{h}\in\mathbb{R}_{+}$
satisfy%
\begin{equation}
||\mathcal{U}^{h};H^{h}||=1,\ ||T^{h}\mathcal{U}^{h}-t^{h}\mathcal{U}%
^{h};H^{h}||:=\delta^{h}\in(0,t^{h}). \label{E6}%
\end{equation}
Then there exists an eigenvalue $\tau_{n}^{h}$ of the operator $T^{h}$ such
that%
\begin{equation}
|t^{h}-\tau_{n}^{h}|\leq\delta^{h}. \label{E7}%
\end{equation}
Moreover, for any $\delta_{\ast}^{h}\in(\delta^{h},t^{h})$, one finds
coefficients $c_{k}^{h}$ to fulfil the relations%
\begin{equation}
\left\Vert \mathcal{U}^{h}-\sum\nolimits_{k=N^{h}}^{N^{h}+X^{h}-1}c_{k}%
^{h}u_{(k)}^{h};H^{h}\right\Vert \leq2\frac{\delta^{h}}{\delta_{\ast}^{h}%
},\ \ \ \ \sum\nolimits_{k=N^{h}}^{N^{h}+X^{h}-1}|c_{k}^{h}|^{2}=1 \label{E8}%
\end{equation}
where $\tau_{N^{h}}^{h},...,\tau_{N^{h}+X^{h}-1}^{h}$are all eigenvalues in
the segment $[t^{h}-\delta_{\ast}^{h},t^{h}+\delta_{\ast}^{h}]$ and
$u_{(N^{h})}^{h},...,u_{(N^{h}+X^{h}-1)}^{h}$ are the corresponding
eigenvectors subject to the normalization and orthogonality conditions%
\begin{equation}
(u_{(p)}^{h},u_{(q)}^{h})_{h}=\delta_{p,q}. \label{E9}%
\end{equation}

\end{lemma}

\subsection{Detecting eigenvalues with prescribed asymptotic
form\label{sect4.4}}

Let $\mu_{p}^{h}\in\mathbb{R}_{+}$ and $v_{(p)}^{h}=(v_{(p)0}^{h},v_{(p)1}%
^{h},...,v_{(p)J}^{h})\in\mathfrak{H}_{-}$ be an eigenvalue of the equation
(\ref{S2}) and the corresponding eigenvector enjoing the normalization and
orthogonality conditions%
\begin{equation}
(v_{(p)}^{h},v_{(q)}^{h})_{\mathfrak{L}}=\delta_{p,q} \label{E10}%
\end{equation}
where $(\ ,\ )_{\mathfrak{L}}$ denotes the scalar product in the Lebesgue
space $\mathfrak{L}$ induced by the norm (\ref{LL}). In Lemma \ref{lemVL} we
set%
\begin{equation}
t_{p}^{h}=(1+\mu_{p}^{h})^{-1} \label{E11}%
\end{equation}
and build an asymptotic approximation $\mathcal{U}_{(p)}^{h}$ of an
eigenfunction of the problem (\ref{7})-(\ref{13}) in the junction (\ref{4}).
To mimic the method of matched asymptotic expansions, we use asymptotic
structures with "overlapping" cut-off functions, see \cite[Ch.2]{MaNaPl},
\cite{na239, na344} etc., namely in addition to the functions $\chi_{j}$ in
(\ref{A5}) and $X_{j}$ in (\ref{S6}) we introduce%
\begin{align}
\mathcal{X}^{h}(y)  &  =0\text{ for }r_{j}\leq
Rh,\ j=1,...,J,\ \ \ \ \mathcal{X}^{h}(y)=1\text{ for }\min\{r_{1}%
,...,r_{J}\}\geq2Rh,\label{E12}\\
X^{h}(z)  &  =0\text{ for }z\leq2h,\ \ \ \ X^{h}(z)=1\text{ for }%
z\geq3h.\nonumber
\end{align}
Radius $R$ is chosen such that $\mathcal{X}^{h}(y)=0$ on $\omega_{j}^{h}.$

The functions $\mathcal{U}_{(p)}^{h}$ and $\mathcal{V}_{(p)}^{h}$ are
determined by formulas%
\begin{align}
\mathcal{U}_{(p)}^{h}  &  =||\mathcal{V}_{(p)}^{h};H^{h}||^{-1}\mathcal{V}%
_{(p)}^{h},\label{E13}\\
\mathcal{V}_{(p)0}^{h}(x)  &  =\mathcal{X}^{h}(y)v_{(p)0}^{h}(y)+%
%TCIMACRO{\tsum \nolimits_{j}}%
%BeginExpansion
{\textstyle\sum\nolimits_{j}}
%EndExpansion
\chi_{j}(y)(c_{j}^{0}+c_{j}^{1}(\mathbf{W}_{j}^{0}(\xi^{j})+h\widehat
{\mathbf{w}}_{j}(\xi^{j}))\label{E14}\\
&  -%
%TCIMACRO{\tsum \nolimits_{j}}%
%BeginExpansion
{\textstyle\sum\nolimits_{j}}
%EndExpansion
\mathcal{X}^{h}(y)\chi_{j}(y)\left(  c_{j}^{0}+c_{j}^{1}\frac{1}{2\pi}%
(\ln\frac{1}{|\eta^{j}|}+\ln c_{\log}(\omega_{j})\right)  ,\nonumber\\
\mathcal{V}_{(p)j}^{h}(x)  &  =X^{h}(z)v_{(p)j}^{h}(z)+X_{j}(z)(c_{j}%
^{0}+hc_{j}^{1}\mathbf{W}_{j}(\xi^{j}))-X^{h}(z)X_{j}(z)(c_{j}^{0}+hc_{j}%
^{1}\gamma_{j}^{-1}|\omega_{j}|^{-1}\zeta^{j}) \label{E15}%
\end{align}
where ingredients are taken from (\ref{B9}), (\ref{B10}) and $\widehat
{\mathbf{w}}_{j}\in H^{1}(\Lambda_{j})$ is a function with compact support
such that%
\begin{equation}
\widehat{\mathbf{w}}_{j}(\xi^{j})=\mathbf{W}_{j}(\xi^{j}),\ \ \xi^{j}%
\in\partial\omega_{j}\times(0,1). \label{wwj}%
\end{equation}
The latter condition and the cut-off functions $\mathcal{X}^{h},X^{h}$ in
(\ref{E14}), (\ref{E15}) assure that $\mathcal{V}_{(p)0}^{h}$ and
$\mathcal{V}_{(p)j}^{h}$ coincide with each other on $\upsilon_{j}(h)$, cf.
(\ref{12}), and the composite function $\mathcal{V}_{(p)}^{h}$ falls into
$H_{0}^{1}(\Xi(h);\Gamma(h)).$

First of all, we compute the scalar products $(\mathcal{V}_{(p)}%
^{h},\mathcal{V}_{(q)}^{h})_{h}$. To this end, we observe that, according to
(\ref{B11}), we have
\begin{equation}%
%TCIMACRO{\tsum \nolimits_{j}}%
%BeginExpansion
{\textstyle\sum\nolimits_{j}}
%EndExpansion
(|b_{(p)j}^{h}|+|\partial_{z}v_{(p)j}^{h}(0)|)\leq c|\ln h|^{-1}%
%TCIMACRO{\tsum \nolimits_{j}}%
%BeginExpansion
{\textstyle\sum\nolimits_{j}}
%EndExpansion
(|\widehat{v}_{(p)0}^{h}(P^{j})|+|v_{(p)j}^{h}(0)|). \label{E16}%
\end{equation}
The estimate (\ref{10}) applied in the problem (\ref{20}), (\ref{21}) shows
that%
\begin{equation}
||\widehat{v}_{(p)0}^{h};H^{1}(\omega_{0})||+%
%TCIMACRO{\tsum \nolimits_{j}}%
%BeginExpansion
{\textstyle\sum\nolimits_{j}}
%EndExpansion
|\widehat{v}_{(p)0}^{h}(P^{j})|\leq c\mu_{p}||v_{(p)0}^{h};L^{2}(\omega
_{0})||\leq C_{p}. \label{E17}%
\end{equation}
Moreover, a solution of (\ref{21}), (\ref{23}) satisfies%
\begin{equation}
||v_{(p)j};H^{2}(I_{j})||\leq c(\mu_{p}||v_{(p)j};L^{2}(I_{j})||+|\partial
_{z}v_{(p)j}(0)|) \label{E18}%
\end{equation}
while the last bound is due to the estimate (\ref{E17}) and the small factor
$|\ln h|^{-1}$ on the right of (\ref{E16}).

Recalling that $c_{j}^{0}=\widehat{v}_{(p)0}^{h}(P^{j})$ and $c_{j}^{1}%
=\gamma_{j}^{-1}|\omega_{j}|\partial_{z}v_{(p)j}(0)$, see (\ref{BBB}), we
rewrite (\ref{E14}) as follows:%
\[
\mathcal{V}_{(p)0}^{h}=\widehat{v}_{(p)0}^{h}+(\mathcal{X}^{h}-1)%
%TCIMACRO{\tsum \nolimits_{j}}%
%BeginExpansion
{\textstyle\sum\nolimits_{j}}
%EndExpansion
\chi_{j}(\widehat{v}_{(p)0}^{h}-\widehat{v}_{(p)0}^{h}(P^{j}))-(2\pi
)^{-1}\mathcal{X}^{h}%
%TCIMACRO{\tsum \nolimits_{j}}%
%BeginExpansion
{\textstyle\sum\nolimits_{j}}
%EndExpansion
\chi_{j}b_{(p)j}^{h}\ln r_{j}+%
%TCIMACRO{\tsum \nolimits_{j}}%
%BeginExpansion
{\textstyle\sum\nolimits_{j}}
%EndExpansion
\chi_{j}c_{j}^{1}(h\widehat{\mathbf{w}}_{j}+\widetilde{\mathbf{W}}_{j}^{0}).
\]
We list the estimates%
\begin{align}
||(1-\mathcal{X}^{h})\chi_{j}(\widehat{v}_{(p)0}^{h}  &  -\widehat{v}%
_{(p)0}^{h}(P^{j}));H^{1}(\Omega_{\bullet}(h))||\leq ch^{3/2-\delta
}||\widetilde{v}_{(p)0}^{h};V_{\delta}^{2}(\omega_{0})||\leq ch,\label{E19}\\
||b_{(p)j}^{h}\chi_{j}\ln r_{j};H^{1}(\Omega_{\bullet}(h))||  &  \leq
ch^{1/2}|\ln h|^{-1}\left(  \int_{h}^{R}(r_{j}^{-2}+1)r_{j}dr_{j}\right)
^{1/2}\leq ch^{1/2}|\ln h|^{-1/2},\nonumber\\
h||\chi_{j}c_{j}^{1}h\widehat{\mathbf{w}}_{j};H^{1}(\Omega_{\bullet}(h))||  &
\leq ch^{3/2}|\ln h|^{-1}||\widehat{\mathbf{w}}_{j};H^{1}(\Lambda
_{j})||,\nonumber\\
||c_{j}^{1}\chi_{j}\widetilde{\mathbf{W}}_{j}^{0};H^{1}(\Omega_{\bullet
}(h))||  &  \leq ch^{1/2}|\ln h|^{-1}\left(  \int_{h}^{R}\left(  \frac{h^{2}%
}{r_{j}^{4}}+\frac{1}{r_{j}^{2}}\right)  r_{j}dr_{j}\right)  ^{1/2}\leq
ch^{1/2}|\ln h|^{-1/2}.\nonumber
\end{align}
These must be commented. Notice that the factor $h^{1/2}$ comes due to
integration in $z\in(0,h)$. In the first estimate we took into account that
$1-\mathcal{X}^{h}=0$ outside the disk $\mathbb{B}_{2Rh}(P^{j})$ where
$r_{j}\leq2Rh$ and applied the weighted inequality (\ref{10}) for
$\widetilde{v}_{(p)0}^{h}(y)=\widehat{v}_{(p)0}^{h}(y)-\widehat{v}_{(p)0}%
^{h}(P^{j})$ with $\delta\in(0,1/2)$. The second and fourth estimates were
derived by a direct calculation of norms and using the decomposition
(\ref{B4}) of $\mathbf{W}_{j}^{0}$ and the bound $c|\ln h|^{-1}$ for
$|c_{j}^{1}|$, cf. (\ref{E16})-(\ref{E18}). The third estimate is obtained by
the coordinate change $x\mapsto\xi^{j}$, see (\ref{B1}).

The above listed estimates support the following relation:%
\begin{gather}
|(\nabla_{x}\mathcal{V}_{(p)0}^{h},\nabla_{x}\mathcal{V}_{(q)0}^{h}%
)_{\Omega_{\bullet}(h)}+(\mathcal{V}_{(p)0}^{h},\mathcal{V}_{(q)0}%
^{h})_{\Omega_{\bullet}(h)}-h(\nabla_{y}\widehat{v}_{(p)0}^{h},\nabla
_{y}\widehat{v}_{(q)0}^{h})_{\omega_{\bullet}(h)}\label{E20}\\
-h(\widehat{v}_{(p)0}^{h},\widehat{v}_{(q)0}^{h})_{\omega_{\bullet}(h)}|\leq
c|\ln h|^{-1/2}.\nonumber
\end{gather}
Moreover,%
\begin{align}
|(\widehat{v}_{(p)0}^{h},\widehat{v}_{(q)0}^{h})_{\omega_{\bullet}%
(h)}-(v_{(p)0}^{h},v_{(q)0}^{h})_{\omega_{0}}|  &  \leq c|\ln h|^{-1/2}%
,\label{E21}\\
|(\nabla_{y}\widehat{v}_{(p)0}^{h},\nabla_{y}\widehat{v}_{(q)0}^{h}%
)_{\omega_{\bullet}(h)}-\mu_{p}^{h}(v_{(p)0}^{h},v_{(q)0}^{h})_{\omega_{0}}|
&  \leq c|\ln h|^{-1/2}. \label{E22}%
\end{align}
Indeed, in (\ref{E21}) we got rid of $\ln r_{j}$ by using the second estimate
(\ref{E19}) and evaluate $||\widehat{v}_{(p)0}^{h};L^{2}(\omega_{j}^{h})||$ by
means of the weighted inequality (\ref{10}) again. To conclude (\ref{E22}), we
observed additionally that%
\[
|(\nabla_{y}\widehat{v}_{(p)0}^{h},\nabla_{y}\widehat{v}_{(q)0}^{h}%
)_{\omega_{0}}-\mu_{p}^{h}(\widehat{v}_{(p)0}^{h},\widehat{v}_{(q)0}%
^{h})_{\omega_{0}}|\leq c|\ln h|^{-1/2}%
\]
because $\widehat{v}_{(p)0}^{h}$ is a solution of the problem%
\[
-\Delta\widehat{v}_{(p)0}^{h}-\mu_{p}^{h}\widehat{v}_{(p)0}^{h}=\widehat
{f}_{(p)0}^{h}\text{ in }\omega_{0},\ \ \ \partial_{\nu}\widehat{v}_{(p)0}%
^{h}=0\text{ on }\partial\omega_{0}%
\]
with a right-hand side which is caused by abolition of $-b_{(p)j}^{h}\chi
_{j}(2\pi)^{-1}\ln r_{j}$ and therefore has the $L^{2}(\omega_{0})$-norm of
order $|\ln h|^{-1}$.

From (\ref{E20})-(\ref{E22}) we derive that%
\begin{equation}
|(\nabla_{x}\mathcal{V}_{(p)0}^{h},\nabla_{x}\mathcal{V}_{(q)0}^{h}%
)_{\Omega_{\bullet}(h)}+(\mathcal{V}_{(p)0}^{h},\mathcal{V}_{(q)0}%
^{h})_{\Omega_{\bullet}(h)}-h(1+\mu_{p}^{h})(v_{(p)0}^{h},v_{(q)0}%
^{h})_{\omega_{0}}|\leq c_{pq}|\ln h|^{-1/2}. \label{E23}%
\end{equation}

In a similar way but with much simpler calculations (recall that $v_{j}^{h}$
is a smooth function on $[0,l_{j}]$) we derive the inequalities%
\begin{equation}
\left\vert h^{-1}\gamma_{j}(\nabla_{x}\mathcal{V}_{(p)j}^{h},\nabla
_{x}\mathcal{V}_{(q)j}^{h})_{\Omega_{j}(h)}+h^{-1}\rho_{j}(\mathcal{V}%
_{(p)j}^{h},\mathcal{V}_{(q)j}^{h})_{\Omega_{j}(h)}-h(1+\mu_{p}^{h}%
)\rho|\omega_{j}|(v_{(p)j}^{h},v_{(q)j}^{h})_{I_{j}}\right\vert \leq c|\ln
h|^{-1/2}. \label{E24}%
\end{equation}
Notice that the factor $h^{-1}$ is compensated due to the relation
$|\omega_{j}^{h}|=h^{2}|\omega_{j}|$. According to (\ref{17}), (\ref{E1}) and
(\ref{E10}), (\ref{LL}), the inequalities (\ref{E23}) and (\ref{E24}) lead us
to%
\begin{equation}
|(\mathcal{V}_{(p)}^{h},\mathcal{V}_{(q)}^{h})_{h}-h(1+\mu_{p}^{h}%
)\delta_{p,q}|\leq c|\ln h|^{-1/2}. \label{E25}%
\end{equation}
Now we evaluate the norm in (\ref{E6}), namely%
\begin{equation}
\delta_{p}^{h}=||T^{h}\mathcal{U}_{(p)}^{h}-t_{p}^{h}\mathcal{U}_{(p)}%
^{h};H^{h}||=(1+\mu_{p}^{h})||\mathcal{V}_{(p)}^{h};H^{h}||^{-1}%
||\mathcal{V}_{(p)}^{h}-(1+\mu_{p}^{h})T^{h}\mathcal{V}_{(p)}^{h};H^{h}||.
\label{dhp}%
\end{equation}
The first factor on the right is bounded, see Theorem \ref{th2T}, and the
second one does not exceed $c_{p}h^{-1}$ owing to (\ref{E25}) with $p=q$.
Using a definition of a Hilbert norm together with formulas (\ref{E1}) and
(\ref{E2}), we obtain that the last factor is equal to%
\begin{equation}
\sup\left\vert (\mathcal{V}_{(p)}^{h}-(1+\mu_{p}^{h})T^{h}\mathcal{V}%
_{(p)}^{h},W^{h})_{h}\right\vert =\sup\left\vert a(\mathcal{V}_{(p)}^{h}%
,W^{h};\Xi(h))-\mu_{p}^{h}b(\mathcal{V}_{(p)}^{h},W^{h};\Xi(h))\right\vert
\label{E26}%
\end{equation}
where supremum is computed over the unit ball in $H^{h}$. Inequality
(\ref{R8}), see \cite[Thm. 9]{BuCaNa1}, indicates bounds for weighted Lebesgue
norms of $\mathcal{W}^{h}$.

We insert representations (\ref{E14}) and (\ref{E15}) into the last
expressions $\mathcal{J}^{h}$ between the modulo sign in (\ref{E26}) and
detach the elementary term%
\begin{equation}
\mathcal{J}_{\mathbf{w}}^{h}=h%
%TCIMACRO{\tsum \nolimits_{j}}%
%BeginExpansion
{\textstyle\sum\nolimits_{j}}
%EndExpansion
c_{j}^{1}(\nabla_{x}\widehat{\mathbf{w}}_{j},\nabla_{x}W^{h})_{\Omega
_{\bullet}(h)}-\mu_{p}^{h}(\widehat{\mathbf{w}}_{j},W^{h})_{\Omega_{\bullet
}(h)},\ \ \ |\mathcal{J}_{\mathbf{w}}^{h}|\leq ch^{3/2}|\ln h|^{-1}.
\label{E27}%
\end{equation}
We derived the estimate in the same way as above, that is, the information on
$\widehat{\mathbf{w}}_{j},$ $c_{j}^{1}$ and the coordinate change $x\mapsto
\xi^{j}$.

Other ingredients are sufficiently smooth while integration by parts and
commuting the Laplace operator with cut-off functions yield%
\begin{align}
\mathcal{J}_{0}^{h}  &  =(\Delta_{x}v_{(p)0}^{h}+\mu_{p}^{h}v_{(p)0}%
^{h},\mathcal{X}^{h}W^{h})_{\Omega_{\bullet}(h)}+([\Delta_{x},\mathcal{X}%
^{h}]\widetilde{v}_{(p)0}^{h},W_{0}^{h})_{\Omega_{\bullet}(h)},\label{E28}\\
\mathcal{J}_{j0}^{h}  &  =(\Delta_{x}\mathbf{W}_{j}^{0},\chi_{j}W_{0}%
^{h})_{\Omega_{\bullet}(h)}+([\Delta_{x},\chi_{j}]\widetilde{\mathbf{W}}%
_{j}^{0},W_{0}^{h})_{\Omega_{\bullet}(h)}\label{E29}\\
&  +c_{j}^{1}\mu_{p}^{h}(\mathbf{W}_{j}^{0}+\mathcal{X}^{h}(2\pi)^{-1}%
(\ln|\eta|-\ln c_{\log}(\omega_{j})),\chi_{j}W^{h})_{\Omega_{\bullet}%
(h)},\nonumber\\
\mathcal{J}_{j}^{h}  &  =h^{-1}(\gamma_{j}\Delta_{x}v_{(p)j}^{h}+\mu_{p}%
^{h}\rho_{j}v_{(p)j}^{h},X^{h}W^{h})_{\Omega_{j}(h)}+h^{-1}\gamma_{j}%
([\Delta_{x},X^{h}]\widetilde{v}_{(p)j}^{h},W^{h})_{\Omega_{j}(h)}%
,\label{E30}\\
\mathcal{J}_{jj}^{h}  &  =c_{j}^{1}\gamma_{j}(\Delta_{x}\mathbf{W}_{j}%
,X_{j}W^{h})+c_{j}^{1}\gamma_{j}([\Delta_{x},X_{j}]\widetilde{\mathbf{W}}%
_{j},W^{h})_{\Omega_{j}(h)}\label{E31}\\
&  +c_{j}^{1}\mu_{p}^{h}\rho_{j}(\mathbf{W}_{j}-X^{h}\gamma_{j}^{-1}%
|\omega_{j}|^{-1}\zeta^{j},W^{h})_{\Omega_{j}(h)}.\nonumber
\end{align}
Note that integrals over the surfaces $\upsilon_{0}(h),$ $\omega_{j}^{0}(0)$
and $\upsilon_{j}(h)$ vanish due to our choice of cut-off functions and
boundary conditions for $v_{0}^{h},$ $v_{j}^{h}$ and $\mathbf{W}_{j}^{0},$
$\mathbf{W}_{j}$, see Section \ref{sect1.2} and \ref{sect3.1}, respectively.
We will estimate all scalar products in (\ref{E28})-(\ref{E30}) and explain
how their sum converts into $\mathcal{J}^{h}-\mathcal{J}_{\mathbf{w}}^{h}$.

In view of the equation (\ref{20}) the first term in (\ref{E28}) is null. The
next term $\mathcal{J}_{0}^{h2}$ in (\ref{E28}) is obtained from the first and
third terms in (\ref{E14}) after commuting the Laplace operator with the
cut-off function $\mathcal{X}^{h}$; notice that%
\begin{align}
\lbrack\Delta_{x},\mathcal{X}^{h}\chi_{j}]  &  =[\Delta_{x},\mathcal{X}%
^{h}]+[\Delta_{x},\chi_{j}],\label{E32}\\
\lbrack\Delta_{x},\mathcal{X}^{h}]  &  =2\nabla_{x}\mathcal{X}^{h}\cdot
\nabla_{x}+\Delta_{x}\mathcal{X}^{h},\ \ \ [\Delta_{x},\chi_{j}]=2\nabla
_{x}\chi_{j}\cdot\nabla_{x}+\Delta_{x}\chi_{j}.\nonumber
\end{align}
Since $\widetilde{v}_{(p)0}^{h}(P^{j})=0$, see (\ref{A9}), and supports of
coefficients in the differential operator $[\Delta_{x},\mathcal{X}^{h}]$
belong to the disk $\mathbb{B}_{2Rh}(P^{j})$, see (\ref{E12}), the direct
consequence of the one-dimensional Hardy inequality%
\begin{align}
||r_{j}^{-2}(1+|\ln r_{j}|)^{-1}\widetilde{v}_{(p)0}^{h};L^{2}(\mathbb{B}%
_{2Rh}(P^{j}))||  &  \leq c||r_{j}^{-2}(1+|\ln r_{j}|)^{-1}\nabla
_{y}\widetilde{v}_{(p)0}^{h};L^{2}(\mathbb{B}_{2Rh}(P^{j}))||\label{EHar}\\
&  \leq||\widetilde{v}_{(p)0}^{h};H^{2}(\mathbb{B}_{2Rh}(P^{j}))||\nonumber
\end{align}
provides that%
\begin{gather}
|\mathcal{J}_{0}^{h1}|\leq(h^{-2}h^{2}(1+|\ln h|)+h^{-1}h(1+|\ln
h|))h^{1/2}||\widetilde{v}_{(p)0}^{h};H^{2}(\mathbb{B}_{2Rh}(P^{j}%
))||\label{E33}\\
\times h(1+|\ln h|)||r^{-1}(1+|\ln r|)^{-1}W_{0}^{h};L^{2}(\Omega_{\bullet
}(h))||\leq ch^{3/2}|\ln h|^{3}.\nonumber
\end{gather}
Here, we took into account that $r=r_{j}<2Rh$ for $y\in\mathbb{B}_{2Rh}%
(P^{j})$. The first (long) multiplier in the middle of (\ref{E33}) is written
according to the relation $|\nabla_{y}^{k}\mathcal{X}^{h}(y)|\leq C_{k}h^{-k}%
$, formula for $[\Delta_{x},\mathcal{X}^{h}]$ in (\ref{E32}) and weights in
(\ref{EHar}). The factor $h^{1/2}$ is due to integration in $z\in(0,h)$ and
finally the weights $r^{-1}(1+|\ln r|)^{-1}$ and $(1+|\ln h|)^{-1}$ from
(\ref{R8}) were considered. It should be also mentioned that $h_{0}<1$ and,
therefore, $|\ln h|\geq c>0$ for $h\in(0,h_{0}]$.

Dealing with (\ref{E29}) we observe that the Laplace equation for
$\mathbf{W}_{j}^{0}$ in (\ref{B3}) annuls the first term $\mathcal{J}%
_{j0}^{h1}$. Supports of coefficients of $[\Delta_{x},\mathcal{\chi}_{j}]$,
see (\ref{E32}), are located in the annulus $\{x\in\overline{\Omega_{0}%
(h)}:R_{j}^{-}\leq r_{j}\leq R_{j}^{+}\}$. Hence, using the decay rate in
$|\eta^{j}|=\rho^{j}=h^{-1}r_{j}$ of the remainder $\widetilde{\mathbf{W}}%
_{j}^{0}$, see (\ref{B4}), we derive the following estimates for the second
and third terms in (\ref{E29}):%
\begin{align}
|\mathcal{J}_{0j}^{h2}|  &  \leq c|c_{j}^{1}|h^{1/2}\left(  \int_{R_{j}^{-}%
}^{R_{j}^{+}}\left(  \frac{1}{h^{2}}\frac{1}{(1+\rho_{j})^{4}}+\frac
{1}{(1+\rho_{j})^{2}}\right)  r_{j}dr_{j}\right)  ^{1/2}||W^{h};L^{2}%
(\Omega_{\bullet}(h))||\leq ch^{3/2},\label{E34}\\
|\mathcal{J}_{0j}^{h3}|  &  \leq c|c_{j}^{1}|h^{1/2}\left(  \int_{Rh}%
^{R_{j}^{+}}\frac{r_{j}dr_{j}}{(1+\rho_{j})^{2}}\right)  ^{1/2}||W^{h}%
;L^{2}(\Omega_{\bullet}(h))||\leq ch^{3/2}|\ln h|.\nonumber
\end{align}
Here, we applied (\ref{R8}) again and recalled that $|c_{j}^{1}|\leq c|\ln
h|^{-1}$. It should be mentioned that $\mathcal{J}_{0j}^{h2}$ and
$\mathcal{J}_{0j}^{h3}$, respectively, involve the commutator $[\Delta
_{x},\mathcal{\chi}_{j}]$ and the multiplication operator $\mu_{p}^{h}$ acting
on $c_{j}^{1}\mathbf{W}_{j}^{0}$ and the last term in (\ref{E14}). In this way
the sum $\mathcal{J}_{0}^{h}+\mathcal{J}_{01}^{h}+...+\mathcal{J}_{0J}%
^{h}+\mathcal{J}_{\mathbf{w}}^{h}$ exhibits the whole part of $\mathcal{J}%
^{h}$ generated by (\ref{E14}).

Referring to (\ref{E30}), we see that $\mathcal{J}_{j}^{h1}=0$ in view of
(\ref{21}). The second term $\mathcal{J}_{j}^{h2}$ in (\ref{E30}) meets the
estimate%
\[
|\mathcal{J}_{j}^{h2}|\leq ch^{-1}(h^{-2}h^{2}+h^{-1}h)h^{3/2}||W_{j}%
^{h};L^{2}(\Omega_{j}^{h})||\leq ch^{3/2}.
\]
Here, $h^{-1}$ came from (\ref{E30}), the relations $|\nabla_{x}^{k}%
X^{h}(z)|\leq c_{k}h^{-k}$ and $\widetilde{v}_{(p)j}^{h}(z)=v_{(p)j}%
^{h}(z)-v_{(p)j}^{h}(0)-z\partial_{z}v_{(p)j}^{h}(0)=O(z^{2})$ were used, the
factor $h^{3/2}$ is due to integration over\ $\Omega_{j}^{h}=\{x\in
\overline{\Omega_{j}(h)}:z\leq3h\}\supset$\textrm{supp}$\partial_{z}X^{h}$ and
finally the direct consequence of the Newton-Leibnitz formula%
\begin{equation}
h^{-2}||W_{j}^{h};L^{2}(\Omega_{j}^{h})||^{2}\leq ch^{-1}(||\partial_{z}%
W_{j}^{h};L^{2}(\Omega_{j}^{h})||^{2}+||W_{j}^{h};L^{2}(\Omega_{j}^{h}%
)||^{2})\leq c||W_{j}^{h};H^{h}||^{2}=c \label{NL}%
\end{equation}
together with the inequality (\ref{R8}) were applied.

Similarly to the above considerations, the first term $\mathcal{J}_{jj}^{h1}$
in (\ref{E31}) vanishes, cf. (\ref{B2}), and the other couple of terms can be
estimated as follows:%
\begin{align*}
|\mathcal{J}_{jj}^{h2}|+|\mathcal{J}_{jj}^{h3}|  &  \leq c|c_{j}^{1}|\left(
\int_{\Omega_{j}^{h}}dx||W_{j}^{h};L^{2}(\Omega_{j}^{h})||^{2}+\int
_{\Omega_{j}(h)}e^{-2\delta z/h}dx||W_{j}^{h};L^{2}(\Omega_{j}(h))||^{2}%
\right)  ^{1/2}\\
&  \leq c|\ln h|^{-1}(h^{3}h^{2}+h^{3}h)^{1/2}\leq c|\ln h|^{-1}h^{2}.
\end{align*}
Here, we took into account the exponential decay in (\ref{B6}) together with
formulas (\ref{NL}) and (\ref{R8}).

In the same way as above we detect a proper redistribution of commutators of
$\Delta_{x}$ with $X^{h},$ $X_{j}$ and conclude that $\mathcal{J}_{j}%
^{h}+\mathcal{J}_{jj}^{h}$ equals a part of $\mathcal{J}^{h}$ generated by
(\ref{E15}).

We summarize our calculations and find that the worst bound in estimates
derived for components of (\ref{E26}), occurs in (\ref{E33}). Hence, according
to (\ref{E25}) with $q=p$, we obtain the following inequality for the quantity
(\ref{dhp}):%
\[
\delta_{p}^{h}\leq c_{p}h^{1/2}|\ln h|^{3}.
\]

Lemma \ref{lemVL} gives us eigenvalues $\tau^{n}(h)$ and $\lambda^{n}%
(h)=\tau^{n}(h)^{-1}-1$, see (\ref{E4}), of the problems (\ref{E3}) and
(\ref{7})-(\ref{13}), respectively, such that, by virtue of (\ref{E11}),%
\begin{align}
\tau^{n}(h)  &  \in\lbrack t_{p}^{h}-c_{p}h^{1/2}|\ln h|^{3},t_{p}^{h}%
+c_{p}h^{1/2}|\ln h|^{3}]\label{Z1}\\
&  \Rightarrow\lambda^{n}(h)\in\left[  \mu_{p}^{h}-C_{p}h^{1/2}|\ln h|^{3}%
,\mu_{p}^{h}+C_{p}h^{1/2}|\ln h|^{3}\right]  \ \ \ \ \forall h\in
(0,h_{p})\nonumber
\end{align}
where $c_{p},$ $C_{p}$ and $h_{p}$ are some positive numbers and the index $n$
may depend on $h$.

\subsection{Theorem on asymptotics\label{sect4.5}}

We are in position to conclude with the main assertion of the paper.

\begin{theorem}
\label{thASYM}For any $N\in\mathbb{N}$, there exist positive $h(N)$ and $c(N)$
such that the entries $\lambda^{1}(h),...,\lambda^{N}(h)$ of the eigenvalue
sequence (\ref{18}) of the original problem (\ref{7})-(\ref{15}) (or
(\ref{16}) in the variational formulation) and the first $N$ positive
eigenvalues%
\begin{equation}
0<\mu_{1}^{h}\leq\mu_{2}^{h}\leq...\leq\mu_{N}^{h} \label{Z2}%
\end{equation}
of the equation (\ref{S1}) or (\ref{S2}) are in the relationship%
\begin{equation}
|\lambda^{n}(h)-\mu_{n}^{h}|\leq c(N)h^{1/2}|\ln h|^{3},\ \ n=1,...,N,\ h\in
(0,h(N)). \label{Z3}%
\end{equation}

\end{theorem}

\textbf{Proof.} The result directly stems from (\ref{Z1}) and all the previous
considerations. We still need to utter several important remarks, in order to
complete the proof. First, if $\mu_{p}^{h}$ is eigenvalue of multiplicity
$\varkappa_{p}^{h}\geq1$, then Lemma \ref{lemVL} provides us with
$\varkappa_{p}^{h}$ different eigenvalues $\lambda^{n}(h),...,\lambda
^{n+\varkappa_{p}^{h}-1}(h)$ satisfying (\ref{Z1}) with a bigger constant
$C_{p}$. Indeed, owing to (\ref{E25}), the constructed approximate
eigenfunctions $\mathcal{U}_{(p)}^{h},...,\mathcal{U}_{(p+\varkappa_{p}%
^{h}-1)}^{h}$ are almost orthonormalized, that is%
\[
|(\mathcal{U}_{(k)}^{h},\mathcal{U}_{(q)}^{h})_{h}-\delta_{p,q}|\leq c_{p}|\ln
h|^{-1/2},\ \ k,q=p,...,p+\varkappa_{p}^{h}-1.
\]
Moreover, setting $\delta_{\ast}^{h}=R\max\{\delta_{p}^{h},...,\delta
_{p+\varkappa_{p}^{h}-1}^{h}\}$ in (\ref{E8}), we see that their projection
onto the linear hull $\mathcal{L}(u_{N^{h}}^{h},...,u_{N^{h}+X^{h}-1}^{h})$
are linear independent for a small $h$ and a big $R$ that is possible in the
case $X^{h}\geq\varkappa_{p}^{h}$ only. Thus, changing $C_{p}\mapsto RC_{p}$
in (\ref{Z1}) gives us at least $\varkappa_{p}^{h}$ desired eigenvalues.

Second, since each eigenvalue $\mu_{p}^{h}$ in (\ref{Z2}) has an eigenvalue
$\lambda^{n^{h}(p)}(h)$ in its small neighborhood and $n(p)\neq n(q)$ for
$p\neq q$, we obtain $\lambda^{n}(h)\leq\mu_{n}^{h}+c_{n}h^{1/2}|\ln h|^{3}$
and confirm the assumption (\ref{R1}) so that Theorem \ref{th2T} becomes true.

Finally, we assume without loss of generality that $N$ is fixed such that the
eigenvalues $\mu_{N}^{0}$ and $\mu_{N+1}^{0}$ of the operator $\mathcal{A}%
^{0}$ obey the relation $\mu_{N}^{0}<\mu_{N+1}^{0}$. If it happens that the
index $n^{h}(N)$ of the eigenvalue $\lambda^{n^{h}(N)}$ in the vicinity of
$\mu_{N}^{h}$ is strictly bigger that $N$, then, for an infinitesimal positive
sequence $\{h_{k}\}_{k\in\mathbb{N}}$, we have $\lambda^{N+1}(h_{k})\leq
\mu_{N+1}^{0}-\varepsilon$ with some $\varepsilon>0$. We can apply Theorem
\ref{th2T} and conclude that $\lambda^{N+1}(0)=\lim\lambda^{N+1}(h_{k})$ as
$k\rightarrow+0$ is an eigenvalue in the interval $(0,\mu_{N+1}^{0}%
-\varepsilon]$, while the limits of (\ref{R2}), (\ref{R3}) constructed from
the corresponding eigenfunctions $u^{N+1}(h_{k},\cdot)$ are orthogonal in
$\mathcal{L}$ to the vector eigenfunctions $(v_{0}^{p},v_{1}^{p},...,v_{J}%
^{p})$, $p=1,...,N,$ of the operator $\mathcal{A}^{0}$. Since $\lambda
^{N+1}(0)$ belongs to the spectrum of $\mathcal{A}^{0}$ but $\lambda
^{N+1}(0)\leq\mu_{N}^{0}$, the latter is absurd. $\blacksquare$

\begin{remark}
\label{remNEGA}A proof of Proposition \ref{prop11S} may follow the same scheme
and meets crucial simplifications because, first, the total multiplicity of
the negative spectrum is known a priori and, second, the corresponding
approximate eigenfunctions (\ref{S5}), (\ref{S6}) decay exponentially at a
distance from the points $P^{1},...,P^{J}$.
\end{remark}

\section{Final remarks\label{sect5}}

\subsection{Simplified and rough asymptotics\label{sect5.1}}

Since the point condition (\ref{A26}) contains the big parameter $\ln h$ in
the matrix (\ref{B12}), it is straightforward to write an asymptotics in $|\ln
h|$ for eigenvalues (\ref{18}). In view of the precision estimate (\ref{Z3})
in Theorem \ref{thASYM} it suffices to find a decomposition of eigenvalues in
the model problem (\ref{20})-(\ref{23}), (\ref{A26}), (\ref{A27}), in
particular, of the projections (\ref{A13}) of the corresponding eigenvectors.

Let us demonstrate the simplest ansatz for the first eigenvalue%
\begin{equation}
\lambda^{1}(h)\sim\mu^{1}(h)\sim|\ln h|^{-1}\mu_{1}^{1}+|\ln h|^{-2}\mu
_{2}^{1}+... \label{71}%
\end{equation}
The corresponding eigenvector of the model problem is searched in the
asymptotic form%
\begin{equation}
v^{1}(h,\cdot)\sim v_{(0)}^{1}+|\ln h|^{-1}v_{(1)}^{1}+|\ln h|^{-2}v_{(2)}%
^{1}+... \label{72}%
\end{equation}
while the absence in (\ref{71}) of the term $|\ln h|^{0}\mu_{0}^{1}$ clearly
requires that%
\begin{equation}
v_{(0)}^{1}=(1,0,...,0) \label{73}%
\end{equation}
and therefore%
\begin{equation}
\wp_{+}^{\prime}v_{(0)}^{1}=\varepsilon=(1,...,1)\in\mathbb{R}^{J},\ \ \wp
_{+}^{\prime\prime}v_{(0)}^{1}=\wp_{-}^{\prime}v_{(0)}^{1}=\wp_{-}%
^{\prime\prime}v_{(0)}^{1}=0\in\mathbb{R}^{J}. \label{74}%
\end{equation}
Hence, a problem to determine $v_{(1)}^{1}=(v_{(1)0}^{1},v_{(1)1}%
^{1},...,v_{(1)J}^{1})$ involves the Neumann problem%
\begin{equation}
-\Delta_{y}v_{(1)0}^{1}(y)=\mu_{1}^{1}v_{(1)0}^{1}(y)=\mu_{1}^{1}%
,\ \ y\in\omega_{\odot},\ \ \ \ \partial_{\nu}v_{(1)0}^{1}(y)=0,\ \ y\in
\partial\omega_{0}, \label{75}%
\end{equation}
which is derived by inserting (\ref{71})-(\ref{73}) into the equations
(\ref{20})-(\ref{22}) and extracting terms of order $|\ln h|^{-1}$.
Furthermore, we obtain from the point condition (\ref{A26}) with the matrix
$S^{h}=(2\pi)^{-1}|\ln h|\mathbb{I}+...$ from (\ref{B12}) that%
\begin{equation}
\wp_{+}^{\prime}v_{(0)}^{1}+(2\pi)^{-1}\wp_{-}^{\prime}v_{(1)}^{1}%
=0\Rightarrow\wp_{-}^{\prime}v_{(1)}^{1}=-2\pi\varepsilon. \label{76}%
\end{equation}
The generalized Green formula (\ref{A14}) delivers the compatibility condition
in this problem, namely%
\begin{align*}
\mu_{1}^{1}|\omega_{0}|  &  =-(\Delta_{y}v_{(1)0}^{1},v_{(0)0}^{1}%
)_{\omega_{0}}+(v_{(1)0}^{1},\Delta_{y}v_{(0)0}^{1})_{\omega_{0}}=\left\langle
\wp_{+}v_{(1)}^{1},\wp_{-}v_{(0)}^{1}\right\rangle -\left\langle \wp
_{-}v_{(1)}^{1},\wp_{+}v_{(0)}^{1}\right\rangle \\
&  =-\left\langle \wp_{-}^{\prime}v_{(1)}^{1},\wp_{+}^{\prime}v_{(0)}%
^{1}\right\rangle =2\pi\left\langle \varepsilon,\varepsilon\right\rangle =2\pi
J.
\end{align*}
Note that we obtained in the ansatz (\ref{71}) the first term%
\begin{equation}
\mu_{1}^{1}=2\pi J|\omega_{0}|^{-1} \label{77}%
\end{equation}
which does not contain much information about the rod elements $\Omega
_{1}(h),...,\Omega_{J}(h)$ of the junction $\Xi(h)$, namely only their number
$J$.

In order to construct the second term $\mu_{2}^{1}$, we, first of all, observe
that a solution of the problem (\ref{75}) can be subject to the orthogonality
condition%
\begin{equation}
\int_{\omega_{0}}v_{(1)0}^{1}(y)dy=0 \label{78}%
\end{equation}
and, thus, formulas (\ref{75}), (\ref{76}) and (\ref{A34}), (\ref{A36}) lead
to the representation%
\[
v_{(1)}^{1}(y)=J^{-1}|\omega_{0}|\mu_{1}^{1}(\mathbf{G}^{1}(y)+...+\mathbf{G}%
^{J}(y))=2\pi(\mathbf{G}^{1}(y)+...+\mathbf{G}^{J}(y))
\]
so that%
\begin{equation}
\wp_{+}^{\prime}v_{(1)}^{1}=2\pi\mathcal{G}\varepsilon,\ \ \ \wp_{+}%
^{\prime\prime}v_{(1)}^{1}=-2\pi\mathcal{Q}\varepsilon. \label{79}%
\end{equation}
Owing to (\ref{78}), the compatibility condition in the problem%
\[
-\Delta_{y}v_{(2)0}^{1}(y)=\mu_{2}^{1}v_{(0)0}^{1}(y)+\mu_{1}^{1}v_{(1)0}%
^{1}(y),\ \ y\in\omega_{\odot},\ \ \ \ \partial_{z}v_{(2)0}^{1}(y)=0,\ \ y\in
\partial\omega_{0},
\]
reads%
\[
\mu_{2}^{1}|\omega_{0}|=-(\Delta_{y}v_{(2)0}^{1},v_{(0)0}^{1})_{\omega_{0}%
}=-\left\langle \wp_{-}^{\prime}v_{(2)}^{1},\wp_{+}^{\prime}v_{(0)}%
^{1}\right\rangle
\]
and implies%
\begin{equation}
\mu_{2}^{1}=2\pi|\omega_{0}|^{-1}\left\langle (C_{\log}+2\pi\mathcal{G+}%
2\pi\mathcal{Q})\varepsilon,\varepsilon\right\rangle . \label{80}%
\end{equation}
The diagonal matrices $\mathcal{Q}$ and $C_{\log}$ depend on length $l_{j}$
and the reduced cross-section $\omega_{j}$ of the rod $\Omega_{j}(h),$
$j=1,...,J$, while the matrix $\mathcal{G}$ reflects disposition of the
junction points $P^{1},...,P^{J}$ in the base of the plate $\Omega_{0}(h)$.

The terms (\ref{77}) and (\ref{80}) detached in (\ref{71}) have been computed.
An estimate of the asymptotic remainder $\widetilde{\mu}^{1}(h)$ is a simple
algebraic task and Theorem \ref{thASYM} converts the estimate into%
\begin{equation}
|\lambda^{1}(h)-|\ln h|^{-1}\mu_{1}^{1}-|\ln h|^{-2}\mu_{2}^{1}|\leq c_{1}|\ln
h|^{-3}. \label{81}%
\end{equation}

\subsection{The homogeneous junction\label{sect5.2}}

By setting in (\ref{14}) the restrictions (\ref{82}), the problem
(\ref{7})-(\ref{13}) reduces to the mixed boundary-value problem%
\begin{align}
-\Delta_{x}u(h,x)  &  =\lambda(h)u(h,x),\ \ \ x\in\Xi_{\sqcup}(h),\label{83}\\
u(h,x)  &  =0,\ x\in\Gamma(h),\ \ \ \partial_{\nu}u(h,x)=0,\ x\in\Xi_{\sqcup
}(h)\setminus\Gamma(h),\nonumber
\end{align}
stated in the intact domain $\Xi_{\sqcup}(h)=\Omega_{0}(h)\cup\Omega
_{1}(h)\cup...\cup\Omega_{J}(h)$, cf. (\ref{4}) and (\ref{1}), (\ref{3}). An
asymptotic analysis of the stationary problem (\ref{83}) with a source term
$f(h,\cdot)\in L^{2}(\Xi_{\sqcup}(h))$ instead of $\lambda(h)u(h,\cdot)$ has
been developed in \cite[Sect. 2]{BuCaNa1}. Let us outline certain
peculiarities of asymptotic models in the case (\ref{82}), which are to be
derived along the same scheme as in Sections \ref{sect3} and \ref{sect4}.

As was mentioned in Section \ref{sect1.2}, a specific feature of the
homogeneous junction $\Xi_{\sqcup}(h)$ is that the limit operator decouples,
see (\ref{88}), i.e., instead of the connected skeleton in fig. \ref{f2},a, we
obtain in the limit the disjoint domain $\omega\subset\mathbb{R}^{2}$ and
intervals $I_{1},...,I_{J}\subset\mathbb{R}$ as drawn in fig. \ref{f2},b.

Both the asymptotic models can be applied for the stationary problem of type
(\ref{7})-(\ref{13}) with the source term $f_{p}$ instead of $\lambda u_{p}$
on the right-hand side of the Poisson equations. However, a result in
\cite{BuCaNa1} displays an explicit rational dependence on $|\ln h|$ of the
corresponding solution that furnishes its complete asymptotic form in finite
steps. Moreover, it reduces an evident importance of the models for the
spectral problem (\ref{7})-(\ref{13}) where an argument based on the big
parameter $|\ln h|$ in (\ref{B12}) and (\ref{A26}) detects the holomorphic
dependence on $|\ln h|^{-1}$ which clearly cannot be described in finite steps.

Trying to formulate point conditions connecting the limit problems (\ref{20}),
(\ref{22}) and (\ref{21}), (\ref{23}) in the skeleton $\Xi^{0}$ of the
junction $\Xi(h)$, see fig. \ref{f1} and \ref{f2}, we need a special solution
of the homogeneous Neumann problem for the Laplace equation in the union
$\Upsilon_{j}=\Lambda\cup Q_{j}$ of a layer and a semi-infinite cylinder. As
was shown in \cite{BuCaNa1}, this solution gets the asymptotic behavior%
\begin{align}
W(\xi)  &  =\frac{1}{2\pi}\left(  \ln\frac{1}{|\eta|}+\ln c_{\log}(\Lambda\cup
Q_{j})\right)  +\sum\limits_{i=1,2}K_{i}(\Lambda\cup Q_{j})\frac{\eta_{j}%
}{|\eta|^{2}}+O(|\eta|^{-2}),\ \xi\in\Lambda,\ |\eta|\rightarrow
+\infty,\label{84}\\
W(\xi)  &  =|\omega_{j}|^{-1}\zeta+O(e^{-\delta\zeta}),\ \delta>0,\ \xi\in
Q_{j},\ \zeta\rightarrow+\infty. \label{85}%
\end{align}
It should be noticed that a proper choice of the coordinate origin eliminates
the constants $K_{i}=K_{i}(\Lambda\cup Q_{j})$ because the change $\eta
\mapsto\eta^{\prime}=\eta+K$ with $K\in\mathbb{R}^{2}$ provides the formulas%
\begin{align*}
\ln|\eta|  &  =\ln|\eta^{\prime}-K|=\ln|\eta^{\prime}|-|\eta^{\prime}%
|^{-2}(K\cdot\eta^{\prime})+O(|\eta^{\prime}|^{-3}),\\
|\eta|^{-2}\eta_{i}  &  =|\eta^{\prime}|^{-2}\eta_{i}^{\prime}+O(|\eta
^{\prime}|^{-3}).
\end{align*}
In what follows we assume that the points $P^{j}$ and, therefore, the
coordinates $y^{j}=h\eta^{j}$ are fixed such that $K_{i}(\Lambda\cup
Q_{j})=0,$ $i=1,2,$ in (\ref{84}).

Repeating the matching procedure performed in Section \ref{sect3.2}, we keep
the expansions (\ref{B7}), (\ref{B8}) but, according to (\ref{84}),
(\ref{85}), replace (\ref{B9}), (\ref{B10}) with the following ones:%
\begin{align*}
c_{j}^{0}+c_{j}^{1}W_{j}(\xi)  &  =c_{j}^{1}\frac{1}{2\pi}\ln\frac{1}%
{|\eta^{j}|}+c_{j}^{1}\frac{1}{2\pi}\ln c_{\log}(\Lambda\cup Q_{j}))+c_{j}%
^{0}+...\text{ \ \ in }\Omega_{\bullet}(h),\\
c_{j}^{0}+c_{j}^{1}W_{j}(\xi)  &  =c_{j}^{0}+c_{j}^{1}|\omega_{j}|^{-1}%
\zeta+...\text{ \ \ in }\Omega_{j}(h).
\end{align*}
Thus, we obtain the relations%
\begin{equation}
b_{j}^{h}-h|\omega_{j}|^{-1}\partial_{z}v_{j}^{h}(0)=0,\ \ \ \widehat{v}%
_{0}^{h}(P^{j})-v_{j}^{h}(0)=b_{j}^{h}(2\pi)^{-1}(\ln h+\ln c_{\log}%
(\Lambda\cup Q_{j})) \label{86}%
\end{equation}
which look quite similar to (\ref{B11}) but we have the small factor $h$ on
the derivatives $\partial_{z}v_{j}^{h}(0)$, that is, on the projection
$\wp_{-}^{\prime\prime}v$. To perform the correct limit passage $h\rightarrow
+0$, we recall our previous asymptotic analysis in \cite[Sect. 2]{BuCaNa1} and
make the substitution%
\begin{equation}
v^{h}=(v_{0}^{h},v_{1}^{h},...,v_{J}^{h})\Rightarrow\mathbf{v}^{h}%
=(\mathbf{v}_{0}^{h},\mathbf{v}_{1}^{h},...,\mathbf{v}_{J}^{h})=(v_{0}%
^{h},h^{-1/2}v_{1}^{h},...,h^{-1/2}v_{J}^{h}). \label{87}%
\end{equation}
As a result, we conclude the point conditions%
\begin{equation}
b_{j}=0,\ v_{j}^{0}(0)=0,\ j=1,...,J\ \ \ \ \Leftrightarrow\ \ \ \ \wp
_{-}^{\prime}\mathbf{v}^{0}=0\in\mathbb{R}^{J},\ \wp_{+}^{\prime\prime
}\mathbf{v}^{0}=0\in\mathbb{R}^{J} \label{88}%
\end{equation}
corresponding to the self-adjoint extension described in Remark \ref{remDIRI}.
In other words, the limit spectrum (\ref{T9}) is composed from the spectrum
$\{\kappa_{h}^{0}\}_{n\in\mathbb{N}}$ of the Neumann problem in $\omega_{0}$
and the spectra $\{\kappa_{h}^{j}=\pi^{2}l_{j}^{-2}n^{2}\}_{n\in\mathbb{N}}$
of the Dirichlet problems in $I_{j},$ $j=1,...,J$ . It is worth to mention
that, as was observed in \cite{BuCaNa1}, the substitution (\ref{87}) has a
clear physical reason, namely the energy functional for the problem (\ref{83})
gets an appropriate approximation by the sum of the energy functionals for the
above mentioned limit problems multiplied with the common factor $h$.

Let us construct the first correction term in the asymptotics%
\begin{equation}
\lambda^{1}(h)=0+h\mu_{1}^{1}+\widetilde{\lambda}^{1}(h) \label{89}%
\end{equation}
of the first eigenvalue of the problem (\ref{83}). Recalling an asymptotic
procedure in \cite[Section 2]{BuCaNa1}, we search for the corresponding
eigenfunction in the form%
\begin{align}
u^{h}(h,x)  &  =1+hv_{(1)0}^{1}(x)+...\ \ \ \text{in }\Omega_{\bullet
}(h),\label{90}\\
u^{h}(h,x)  &  =1-l_{j}^{-1}z+...\ \ \ \text{in }\Omega_{j}(h).\nonumber
\end{align}
Regarding (\ref{90}) as outer expansions, we write the inner expansions in the
vicinity of the sockets $\theta_{j}^{h}$ as follows:%
\[
u^{h}(h,x)=1-hl_{j}^{-1}|\omega_{j}|W_{j}(h^{-1}(y-P^{j}),h^{-1}z)+...
\]
Finally, we take the representations (\ref{84}), (\ref{85}) and apply the
matching procedure, cf. Section \ref{sect3.2}, to close the problem (\ref{75})
with the asymptotic conditions near the points $P^{1},...,P^{J}$%
\begin{equation}
v_{(1)0}^{1}(x)=-\frac{1}{2\pi}\sum\nolimits_{j}\chi_{j}(y)\frac{|\omega_{j}%
|}{l_{j}}\ln\frac{1}{r_{j}}+\widehat{v}_{(1)0}^{1}(x),\ \ \ \widehat{v}%
_{(1)0}^{1}\in H^{2}(\omega_{0}). \label{91}%
\end{equation}
Similarly to Section \ref{sect5.1} the compatibility condition in the problem
(\ref{75}), (\ref{91}) converts into the formula%
\[
\mu_{1}^{1}=|\omega_{0}|^{-1}(|\omega_{1}|l_{1}^{-1}+...+|\omega_{J}%
|l_{J}^{-1}).
\]
Combining approaches in Section \ref{sect4} and \cite[Sect. 2]{BuCaNa1}, the
asymptotic formula (\ref{89}) can be justified by means of the estimate
$|\widetilde{\lambda}^{1}(h)|\leq c_{1}h^{2}(1+|\ln h|)$ for the remainder.

Let us consider the problems (\ref{20}), (\ref{22}) and (\ref{21}), (\ref{23})
connected through the point conditions%
\begin{equation}
\wp_{-}^{\prime}\mathbf{v}^{h}+h^{1/2}\wp_{-}^{\prime\prime}\mathbf{v}%
^{h}=0,\ \ \ h^{-1/2}\wp_{+}^{\prime\prime}\mathbf{v}^{h}-\wp_{+}^{\prime
}\mathbf{v}^{h}-S^{h}\wp_{-}^{\prime}\mathbf{v}^{h}=0\in\mathbb{R}^{J},
\label{93}%
\end{equation}
where $\gamma_{j}=\rho_{j}=1$ because the junction is homogeneous and the
$J\times J$-matrix $S^{h}$ is given by formula (\ref{B12}) with $c_{\log
}(\Lambda\cup Q_{j})$ instead of $c_{\log}(\omega_{j})$, see (\ref{84}) and
(\ref{B4}).

The conditions (\ref{93}) follow immediately from (\ref{86}) after
substitution (\ref{87}). They are involved into the symmetric generalized
Green formula of type (\ref{A33})%
\begin{align*}
q(v^{h},w^{h})  &  =\left\langle \wp_{+}^{\prime}v^{h}-h^{-1/2}\wp_{+}%
^{\prime\prime}v^{h}+S^{h}\wp_{-}^{\prime}v^{h},\wp_{-}^{\prime}%
w^{h}\right\rangle -\left\langle \wp_{-}^{\prime}v^{h},\wp_{+}^{\prime}%
w^{h}-h^{-1/2}\wp_{+}^{\prime\prime}w^{h}+S^{h}\wp_{-}^{\prime}w^{h}%
\right\rangle \\
&  +h^{-1/2}\left\langle \wp_{+}^{\prime\prime}v^{h},\wp_{-}^{\prime}%
w^{h}+h^{1/2}\wp_{-}^{\prime\prime}w^{h}\right\rangle -h^{-1/2}\left\langle
\wp_{-}^{\prime}v^{h}+h^{1/2}\wp_{-}^{\prime\prime}v^{h},\wp_{+}^{\prime
\prime}w^{h}\right\rangle ,
\end{align*}
and therefore, all requirements in Section \ref{sect2} and \ref{sect3} with
slight modifications apply to the model in the case (\ref{82}), too. Moreover,
a simple analysis requiring only for algebraic operations as in Section
\ref{sect5.1}, provides the representation%
\begin{equation}
\mu^{1}(h)=h\mu_{1}^{1}+O(h^{2}(1+|\ln h|)) \label{94}%
\end{equation}
of an eigenvalue in the model (\ref{20})-(\ref{23}), (\ref{93}) which is
supplied with an operator of type (\ref{A31}) on the function space
(\ref{A12}) with detached asymptotics.

\bigskip

\textbf{Acknowledgements.}

\bigskip

The work of S.A.N. was supported by grant 15-01-02175 of Russian Foundation
for Basic Research. G.C. is a member of GNAMPA\ of INDAM.


\bigskip

\begin{thebibliography}{99}                                                                                               %


\bibitem {BeFa}Beresin F.A., Faddeev L.D., A remark on Schr\"{o}dinger
equation with a singular potential, \textit{Sov. Math. Doklady}, \textbf{137}
(5) (1961) 1011-1014.

\bibitem {BiSo}Birman M.Sh., Solomjak M.Z., \textit{Spectral theory of
selfadjoint operators in Hilbert space. Mathematics and its Applications}
(Soviet Series). D. Reidel Publishing Co., Dordrecht, 1987.

\bibitem {Pank}Br\"{u}ning J., Geyler V., Pankrashkin K., Spectra of
self-adjoint extensions and applications to solvable Schr\"{o}dinger
operators, \textit{Rev. Math. Phys.} \textbf{20} (1) (2008) 1-70.

\bibitem {BuCaNa1}Bunoiu R., Cardone G., Nazarov S.A., Scalar boundary value
problems on junctions of thin rods and plates. I. Asymptotic analysis and
error estimates, \textit{ESAIM: Math. Modell. Num. Anal.} \textbf{48} (2014) 1495-1528.

\bibitem {ButCaNa1}G. Buttazzo, G.Cardone, S.A.Nazarov, Thin Elastic Plates
Supported over Small Areas. I: Korn's Inequalities and Boundary Layers,
\textit{J. Convex Analysis} \textbf{23} (1) (2016), 347-386.

\bibitem {ButCaNa2}G. Buttazzo, G.Cardone, S.A.Nazarov, Thin Elastic Plates
Supported over Small Areas. II: Variational-asymptotic models, \textit{J.
Convex Analysis} \textbf{24} (3) (2017) 819-855.

\bibitem {Yves}Colin De Verdi\`{e}re Y., Pseudo-Laplaciens II, \textit{Ann.
Inst. Fourier} \textbf{33} (2) (1983) 87-113.

\bibitem {Ko}Kondratiev V.A., Boundary problems for elliptic equations in
domains with conical or angular points, \textit{Trans. Moscow Math. Soc.}
\textbf{16} (1967) 227-313.

\bibitem {KoMaMo}Kozlov V., Maz'ya V., Movchan A., \textit{Asymptotic analysis
of fields in multi-structures.} Oxford Mathematical Monographs. Oxford Science
Publications. The Clarendon Press, Oxford University Press, New York, 1999.

\bibitem {KoMaRo1}Kozlov V.A., Maz'ya V.G., Rossmann J., \textit{Elliptic
boundary value problems in domains with point singularities.} Providence:
Amer. Math. Soc., 1997.

\bibitem {Ilin}Il'in A.M., \textit{Matching of asymptotic expansions of
solutions of boundary value problems,} Translations of Mathematical
Monographs, 102. American Mathematical Society, Providence, RI, 1992.

\bibitem {Land}Landkof N.S., \textit{Foundations of modern potential theory.}
Die Grundlehren der mathematischen Wissenschaften, Band 180. Springer-Verlag,
New York-Heidelberg, 1972.

\bibitem {Leng}Lang S., \textit{Algebra}, Graduate text in Mathematics 211,
Springer 2002.

\bibitem {Lions}Lions J.-L., \textit{Some more remarks on boundary value
problems and junctions. Asymptotic methods for elastic structures }(Lisbon,
1993), 103--118, de Gruyter, Berlin, 1995.

\bibitem {LiMa}Lions J.L., Magenes E., \textit{Non-homogeneous boundary value
problems and applications}, Springer-Verlag, New York-Heidelberg, 1972.

\bibitem {MaNaPl}Maz'ya V.G., Nazarov S.A., Plamenevskij B.A.,
\textit{Asymptotic theory of elliptic boundary value problems in singularly
perturbed domains,} Operator Theory: Advances and Applications, \textbf{112},
Birkh\"{a}user Verlag, Basel (2000).

\bibitem {na285}Nazarov S.A., Junctions of singularly degenerating domains
with different limit dimensions. 2, \textit{Trudy seminar. Petrovskii}.
\textbf{20} (1997) 155-195 (English transl.: \textit{J. Math. Sci.},
\textbf{97} (3) (1999) 155--195.)

\bibitem {na344}Nazarov S.A., Estimates for the accuracy of modeling boundary
value problems on the junction of domains with different limit dimensions,
\textit{Izv. Math.}, \textbf{68} (6) (2004) 1179--1215.

\bibitem {na389}Nazarov S.A. Asymptotic behavior of the solution and the
modeling of the Dirichlet problem in an angular domain with rapidly
oscillating boundary, \textit{St. Petersburg Math. J.}, \textbf{19} (2) (2007) 297-326.

\bibitem {na576}Nazarov S.A., Modeling of a singularly perturbed spectral
problem by means of self-adjoint extensions of the operators of the limit
problems, \textit{Funkt. Anal. i Prilozhen.} \textbf{49} (1) (2015) 31-48
(English transl.: \textit{Funct. Anal. Appl.} \textbf{49} (1) (2015) 25--39).

\bibitem {na345}Nazarov S.A., Elliptic boundary value problems on hybrid
domains, \textit{Funkt. Anal. i Prilozhen}, \textbf{38} (4) (2004), 55-72
(English transl.: \textit{Funct. Anal. Appl.} \textbf{38} (4) (2004) 283-297).

\bibitem {na239}Nazarov S.A., Asymptotic conditions at a point, self-adjoint
extensions of operators and the method of matched asymptotic expansions,
\textit{Trans. Am. Math. Soc.} \textbf{193} (1999) 77-126.

\bibitem {na188}Nazarov S.A. Asymptotic solution to a problem with small
obstacles, \textit{Differential equations.} \textbf{31} (6) (1995) 965-974.

\bibitem {na159}Nazarov S.A., Plamenevskii B.A. Selfadjoint elliptic problems
with radiation conditions on the edges of the boundary, \textit{St. Petersburg
Math. J.} \textbf{4} (3) (1993) 569-594.

\bibitem {NaPl}Nazarov, S.A., Plamenevsky, B.A. \textit{Elliptic problems in
domains with piecewise smooth boundaries.} de Gruyter Expositions in
Mathematics, 13. Walter de Gruyter \& Co., Berlin, 1994.

\bibitem {na165}Nazarov S.A., Plamenevskii B.A. Elliptic problems with
radiation conditions on edges of the boundary, \textit{Sb. Math. }\textbf{77}
(1) (1994) 149-176.

\bibitem {na161}Nazarov S.A., Plamenevskii B.A. A generalized Green's formula
for elliptic problems in domains with edges, \textit{J. Math. Sci.
}\textbf{73} (6) (1995) 674-700.

\bibitem {na504}Nazarov S.A., Specovius-Neugebauer M., Modeling of cracks with
nonlinear effects at the tip zones and the generalized energy criterion,
\textit{Arch. Rational Mech. Anal.} \textbf{202} (2011) 1019-1057.

\bibitem {Pav}Pavlov B.S., The theory of extensions, and explicit solvable
models, \textit{Russian Mathematical Surveys} \textbf{42} (6) (1987) 127-168.

\bibitem {PoSe}P\'{o}lya G., Szeg\"{o} G., \textit{Isoperimetric Inequalities
in Mathematical Physics}, Annals of Mathematics Studies, n. 27, Princeton
University Press, Princeton, N. J., 1951.

\bibitem {Rofe}Rofe-Beketov F.S., Self-adjoint extensions of differential
operators in the space of vector functions,\ \textit{Dokl. Akad. Nauk SSSR}
\textbf{184} (1969) 1034-1037; English transl.: \textit{Soviet Math.Dokl.}
\textbf{10} (1969) 188-192.

\bibitem {Smir}Smirnov V.I., \textit{A course of higher mathematics. }Vol. II.
Advanced calculus. Translation edited by I. N. Sneddon Pergamon Press,
Oxford-Edinburgh-New York-Paris-Frankfurt; Addison-Wesley Publishing Co.,
Inc., Reading, Mass.-London 1964.

\bibitem {VanDyke}Van Dyke M., \textit{Perturbation methods in fluid
mechanics}, Applied Mathematics and Mechanics, Vol. 8 Academic Press, New
York-London 1964.

\bibitem {ViLu}Visik M. I., Ljusternik L.A., Regular degeneration and boundary
layer of linear differential equations with small parameter, \textit{Amer.
Math.Soc. Transl.} \textbf{20} (2) (1962) 239-364.
\end{thebibliography}
\end{document}